\documentclass[11pt,twoside,leqno]{aomamlt2e} 
\input ifpdf.sty
  \pageno{513}
\received{September 5, 2003}
\revised{January 25, 2005}
\newcommand\hfq{\hfill\qed}
\newcommand\nn{\nonumber}
\numberwithin{equation}{section}
  \newcounter{notes}

\newtheorem{theorem}{Theorem}[section]
\newtheorem{proposition}[theorem]{Proposition}
\newtheorem{lemma}[theorem]{Lemma}
\newtheorem{corollary}[theorem]{Corollary}
\newtheorem{conjecture}[theorem]{Conjecture}

\theorembodyfont{\upshape}

\newtheorem{definition}[theorem]{{\it Definition}}

\newcommand{\Ad}{\operatorname{Ad}}

\newcommand{\h}{\textsf{\textup{h}}}
\newcommand{\htop}{\textsf{\textup{h}}_{\mathrm{top}}}
\newcommand{\Hh}{\textsf{\textup{H}}}

\newcommand{\Id}{\operatorname{I}_k}

\newcommand{\GL}{\operatorname{GL}}
\newcommand{\SL}{\operatorname{SL}}
\newcommand{\Mat}{\operatorname{Mat}}

\newcommand{\supp}{\operatorname{supp}}
\newcommand{\dm}{\operatorname{d}\!}

\newcommand{\bs}{\mathbf{s}}
\newcommand{\bt}{\mathbf{t}}

\newcommand{\by}{\mathbf{y}}

\newcommand{\n}{\mathbf{n}}

\newcommand{\HH}{\mathbb{H}}
\newcommand{\RR}{\mathbb{R}}

\newcommand{\QQ}{\mathbb Q}

\newcommand{\ZZ}{\mathbb{Z}}

\newcommand{\cA}{\mathcal{A}}
\newcommand{\cB}{\mathcal{B}}
\newcommand{\cC}{\mathcal{C}}

\newcommand{\cQ}{\mathcal{Q}}
\newcommand{\cP}{\mathcal{P}}

\renewcommand{\cH}{\mathcal{H}}

\newcommand {\absolute}[1] {\left| {#1} \right|}
\newcommand {\norm}[1] {\left\| {#1} \right\|}
\newcommand {\diag}{\operatorname{diag}}
\newcommand {\centralizer}[1] {{C({#1})}}

\newcommand{\IGNORE}[1] {}

\newcommand\Bad{\Xi}

\newsymbol\emptyphi 203F
\renewcommand{\emptyset}{\emptyphi}
\newcommand{\twobytwo}[4] {\begin{pmatrix} {#1}&{#2}\\{#3}&{#4}\end{pmatrix}}

\begin{document}
\currannalsline{164}{2006} 

\title{Invariant measures and the set of\\ exceptions to
Littlewood's conjecture}

 \acknowledgements{A.K. was partially supported by NSF grant   DMS-007133. E.L. was partially 
supported by NSF grants DMS-0140497 and DMS-0434403. Part of the research was conducted while E.L.
was a Clay Mathematics Institute Long Term Prize fellow. Visits of A.K. and E.L. to the University of
Washington were supported by the American Institute of Mathematics and NSF Grant DMS-0222452.}
\twoauthors{Manfred Einsiedler, Anatole Katok,}{Elon Lindenstrauss}

 \institution{University of Washington, Seattle, WA, Princeton University, Princeton, NJ, and Clay Mathematics Institute, Cambridge, MA
\\
\current{Ohio State University, Columbus, OH}\\
\email{manfred@math.osu.edu}
\\ \vglue-9pt
The Pennsylvania State University, State College, PA\\
\email{katok\_a@math.psu.edu}\\
\vglue-9pt
Princeton University, Princeton, NJ\\
\email{elonl@math.princeton.edu}}

 \shorttitle{The set of  exceptions to
Littlewood's conjecture}
 
\centerline{\bf Abstract}
\vglue12pt
We classify the measures on $\SL (k,\Bbb R) / \SL (k, \ZZ)$ which
are invariant and ergodic under the action of the group $A$ of
positive diagonal matrices with positive entropy. We apply this to
prove that the set of exceptions to Littlewood's conjecture has
Hausdorff dimension zero.

\section{Introduction}\label{sec: introduction}

\Subsec{Number theory and dynamics}\label{intro number}
There is a long and rich tradition of applying dynamical methods
to number theory. In many of these applications, a key role is
played by the space $\SL (k, \RR) / \SL (k, \ZZ)$ which can be
identified as the space of unimodular lattices in $\RR ^ k$. Any
subgroup $H < \SL (k, \RR)$ acts on this space in a natural way,
and the dynamical properties of such actions often have deep
number theoretical implications.

A significant landmark in  this direction
is the solution by  G.\ A.\ Margulis \cite{Margulis-Oppenheim-proof}
of the long-standing Oppenheim Conjecture through the study of the
action of a certain subgroup  $H$ on the space of unimodular
lattices in three space. This conjecture, posed by A. Oppenheim in
1929, deals with density properties of the values of indefinite
quadratic forms in three or more variables. So far there is no
proof known of this result in its entirety which avoids the use of
dynamics of homogeneous actions.

An important property of the acting group $H$ in the case of the Oppenheim
Conjecture is that it is generated by unipotents: i.e.\ by elements of $\SL (k,
\RR)$ all of whose eigenvalues are 1. The dynamical result proved by Margulis
was a special case of a conjecture of M.\ S.\ Raghunathan regarding the actions
of general unipotents groups. This  conjecture (and related conjectures made
shortly thereafter) state that for the action of  $H$ generated by unipotents
by left translations on the homogeneous space $G/\Gamma$ of an arbitrary
connected Lie group $G$  by a lattice $\Gamma$,  the only possible $H$-orbit
closures and $H$-ergodic probability measures are of an algebraic type.
Raghunatan's conjecture was proved  in full  generality by M.\ Ratner in a
landmark series of papers (\cite{Ratner-measure-rigidity}, \cite{Ratner91t} and
others; see also the expository papers \cite{Ratner-survey}, \cite{Ratner-SL2}, and
the  
 book \cite{Witte-book}) which led to numerous applications; in
particular, we use Ratner's work heavily in this paper. Ratner's theorems provide
the model for the global orbit structure for systems with {\it parabolic}
behavior. See \cite{HasselblattKatokSurvey} for a general discussion of
principal types of orbit behavior in dynamics.

\Subsec{Weyl chamber  flow and Diophantine approximation}
\label{intro diophantine} In this paper we deal with a different homogeneous
action, which is not so well understood, namely the action by left
multiplication of the group $A$ of positive diagonal $k \times k$ matrices on
$\SL (k, \RR) / \SL (k, \ZZ)$; $A$ is a split Cartan subgroup of   $\SL (k,
\RR)$ and the action of $A$ is also known as a particular case of a {\em  Weyl
chamber flow}~\cite{KatokSpatzierCocycles}.

For $k=2$ the acting group is isomorphic to $\RR$ and the Weyl
chamber flow reduces to the geodesic flow on a surface of constant
negative curvature, namely the modular surface.
This flow  has {\it hyperbolic} structure; it is Anosov
if one makes  minor allowances for  noncompactness and elliptic points.
The orbit structure of such flows is well understood; in particular there is a
great variety of invariant ergodic measures and orbit closures.
For $k > 2$,  the Weyl chamber flow is hyperbolic as an $\RR ^{k-1}$-action,
i.e.\ transversally to the orbits. Such actions are very different from
Anosov flows and display many rigidity properties;
see e.g.\ \cite{KatokSpatzierCocycles}, \cite{KatokSpatzierRigidity}.
One of the manifestations of   rigidity  concerns invariant measures.
Notice that one--parameter subgroups of the Weyl chamber flow are
{\it partially hyperbolic} and each such subgroup
still has many invariant measures.
However, it is conjectured  that $A$-ergodic measures are rare:

\begin{conjecture}[Margulis] \label{special Margulis conjecture}
Let $\mu$ be an $A$-invariant and ergodic probability measure on
$X = \SL (k, \RR) / \SL (k, \ZZ)$ for $k \geq 3$. Then $ \mu$ is
algebraic\/{\rm ;} i.e.\ there is a closed{\rm ,} connected group $L > A$ so
that $\mu$ is the $L$-invariant measure on a single{\rm ,} closed $
L$-orbit.
\end{conjecture}

This conjecture is a special case of much more general conjectures
in this direction by Margulis \cite{Margulis-conjectures}, and by
A.~Katok and R.~Spatzier \cite{KatokSpatzier96}. This type of
behavior was first observed by Furstenberg \cite{Furstenberg-1967}
for the action of the multiplicative semigroup $ \Sigma _ {m,n}= \left\{ m ^ kn ^ l \right\} _ {k,l \geq ^ 1}$
on $\RR / \ZZ$,
where $n,m$ are two multiplicatively independent integers (i.e. not powers of the same integer), and the action is given by $k.x = kx \bmod 1$ for any $k \in \Sigma _ {m, n}$ and $x \in \RR / \ZZ$. Under these assumptions Furstenberg
proved that the only infinite closed invariant set under the
action of this semigroup is the space $ \RR / \ZZ$ itself. He also
raised the question of extensions, in particular to the measure
theoretic analog as well as to the locally homogeneous context.

There is an intrinsic difference regarding the classification of invariant
measures between Weyl chamber flows (e.g.\ higher rank Cartan actions) and
unipotent actions. For unipotent actions, every element of the action already
acts in a rigid manner. For Cartan actions, there is no rigidity for the action
of individual elements, but only for the full action. In stark contrast to
unipotent actions, M.\ Rees \cite{Rees}, \cite[\S9]{EinsiedlerKatok02}
has shown there are lattices $\Gamma < \SL (k, \RR)$ for which there are
nonalgebraic $A$-invariant and ergodic probability measures on $X= \SL (k,
\RR) / \Gamma$ (fortunately, this does not happen for $\Gamma = \SL (k, \ZZ)$,
see \cite{Lindenstrauss-Weiss}, \cite{Margulis-conjectures} and more generally
\cite{Tomanov-maxtori} for related results). These nonalgebraic measures arise precisely because
one-parameter subactions are not rigid, and come from $A$ invariant homogeneous
subspaces which have algebraic factors on which the action degenerates to a one-parameter action.
\IGNORE{these examples come from $A$ invariant homogeneous subspaces of $X$ for which there is an
algebraic projection which intertwines the $A$-action with the action of a one-parameter group (or more
precisely, with a $\RR ^ {k - 1}$ action which is effectively an action of $\RR$)}

While Conjecture~\ref{special Margulis conjecture} is a special case  of
the general question about the structure of invariant measures for  higher rank
hyperbolic  homogeneous actions,  it is of particular interest in view of
number theoretic consequences.
In particular,  it implies the following well-known and long-standing conjecture
of Littlewood \cite[\S2]{Margulis-97}:

\begin{conjecture} [Littlewood (c.\ 1930)]
For every $u,v \in \RR${\rm ,}
\begin{equation} \label{eq: Littlewood}
\liminf _ {n \to \infty} n \langle {n u} \rangle \langle{n
v}\rangle = 0,\end{equation} where $\langle
w \rangle=\min_{n \in \ZZ}|w-n|$ is the distance of $w \in \RR$ to the
nearest integer.
\end{conjecture}

In this paper we prove the following partial result towards
Conjecture~\ref{special Margulis conjecture} which has implications toward
Littlewood's conjecture:

\begin{theorem} \label{theorem about lattice space}
Let $\mu$ be an $A$-invariant and ergodic measure on $X = \SL (k,
\RR) / \SL (k, \ZZ)$ for $k \geq 3$. Assume that there is some one-parameter subgroup of $A$ which acts on
$X$ with positive entropy. Then $\mu$ is algebraic.
\end{theorem}

In \cite{Lindenstrauss-Weiss} a complete classification of the
possible algebraic $\mu$ is given. In particular, we have the following:

\begin{corollary} \label{corollary about lattice space}
Let $\mu$ be as in Theorem~{\rm \ref{theorem about lattice space}.}
Then $\mu$ is not compactly supported.
Furthermore{\rm ,} if $k$ is prime{\rm ,} $\mu$ is the unique $\SL (k, \RR)$-invariant measure on $X$.
\end{corollary}

Theorem~\ref{theorem about lattice space} and its corollary have
the following implication toward\break Littlewood's conjecture:

\begin{theorem}
\label{thm: zero Hausdorff} Let
$$ \Bad = \left\{ (u,v) \in \RR ^ 2:
\liminf _ {n \to \infty} n \langle {n u} \rangle \langle{n
v}\rangle > 0 \right\}
.$$ Then the Hausdorff dimension $\dim_H
\Bad=0$. In fact{\rm ,} $\Bad$ is a countable union of compact sets with
box dimension zero.
\end{theorem}

J.\ W.\ S.\ Cassels and H.\ P.\ F.\ Swinnerton-Dyer
\cite{Cassels-Swinnerton-Dyer} showed that \eqref{eq: Littlewood}
holds for any $u, v$ which are from the same cubic number field
(i.e.\ any field $K$ with degree $[K: \QQ] = 3 $).

It is easy to see that for a.e.\ $(u,v)$ equation \eqref{eq:
Littlewood} holds --- indeed, for almost every $u$ it is already
true that $ \liminf _ {n \to \infty} n \langle {n u} \rangle = 0
$. However, there is a set of $u$ of Hausdorff dimension 1 for
which $ \liminf _ {n \to \infty} n \langle {n u} \rangle>0$; such
$u$ are said to be badly approximable. Pollington and Velani
\cite{Pollington-Velani-00} showed that for every $u \in \RR$, the intersection of the
set
\begin{equation}\label{e: PV}
\{ v \in \RR:(u,v)\mbox{ satisfies
\eqref{eq: Littlewood}}\}
\end{equation}
with the set of badly approximable numbers has Hausdorff dimension one.
Note that this fact is an immediate corollary of
our Theorem \ref{thm: zero Hausdorff} --- indeed, Theorem~\ref{thm: zero
Hausdorff} implies in particular that the complement of this set \eqref{e: PV}
has Hausdorff dimension zero for all $u$. 
We remark that the proof of Pollington and Velani is effective.

Littlewood's conjecture is a special case of a more general question.
More generally, for any $k$ linear forms
$m_i(x _ 1, x _ 2, \dots, x _ k) = \sum_ {j = 1 } ^ k m _ {ij} x _ j$,
one may consider the product
$$
f_m (x _ 1, x _ 2, \dots, x _ k) = \prod_ {i=1}^ k m _ i (x _ 1,
\dots, x _ k),
$$
where $m=(m_{ij})$ denotes the $k \times k$ matrix whose rows are
the linear forms above. Using Theorem~\ref{theorem about lattice
space} we prove the following:
\begin{theorem}\label{theorem about forms}
There is a set $\Xi _ k \subset \SL (k, \RR)$ of Hausdorff dimension $k-1$
so that for every $m \in \SL (k, \RR) \setminus \Xi _ k,$
\begin{equation}
\label{e: product of forms}
\inf _ {\mathbf x \in \ZZ ^ k \setminus \{{\mathbf 0}\}}
|f _ m (\mathbf x)| = 0 .\end{equation} Indeed{\rm ,} this set $\Xi _ k$
is $A$-invariant{\rm ,} and has zero Hausdorff dimension transversally
to the $A$-orbits.
\end{theorem}

For more details, see Section  \ref{sec: transversal} and Section \ref{sec: Littlewood}.
Note that
\eqref{e: product of forms} is automatically satisfied if zero is
attained by $f_m$ evaluated on $\ZZ ^ k \setminus \left\{ 0
\right\}$.

We also want to mention another application of our results due to Hee Oh \cite{Hee-unpublished}, which is related to the following conjecture of Margulis:

\begin{conjecture}[Margulis, 1993] \label{Margulis other conjecture}
Let $G$ be the product of $n \geq 2$ copies of $\SL (2, \RR)${\rm ,}
 \begin{eqnarray*}
U_1& =& \left\{ \twobytwo 1 {*} 0 1 \times \dots \times \twobytwo 1 {*} 0 1 \right\}\\
\noalign{\noindent  and}
 U_2 &=&
\left\{
\twobytwo 1 0 {*} 1 \times \dots \times \twobytwo 1 0 {*} 1 \right\}.
\end{eqnarray*}
 Let $\Gamma < G$ be a discrete
subgroup so that for both $i = 1$ and $2${\rm ,}
 $\Gamma \cap U _ i$ is a lattice in $U _ i$ and for any proper
connected normal subgroup $N<G$ the intersection
 $\Gamma \cap N \cap U _ i$ is trivial. Then $\Gamma$
is commensurable with a Hilbert modular lattice\footnote{For a definition of Hilbert modular lattices, see
\cite{Hee-TIFR}.} up to conjunction in $\GL (2, \RR) \times \dots \times \GL (2, \RR)$.
\end{conjecture}

Hee Oh \cite{Hee-TIFR}  has shown that assuming a topological analog to Conjecture~\ref{special Margulis conjecture} (which is implied by  Conjecture~\ref{special Margulis conjecture}), Conjecture~\ref{Margulis other conjecture} is true for $n \geq 3$. As explained in
\cite{Hee-unpublished} (and following directly from \cite [Thm. 1.5]{ Hee-TIFR}), our result,
Theorem~\ref{theorem about lattice space}, implies the following weaker result (also for $n \geq 3$):
consider the set $\mathcal{D}$ of possible intersections $\Gamma \cap U _ 1$ for $\Gamma$ as in
Conjecture~\ref{Margulis other conjecture}, which is a subset of the space of lattices in $U _ 1$. This set
$\mathcal{D}$ is clearly invariant under conjugation by the diagonal group in $\GL (2, \RR) \times \dots
\times \GL (2, \RR)$; Theorem~\ref{theorem about lattice space} (or more precisely Theorem~\ref{thm:
transversal-2} which we prove using Theorem~\ref{theorem about lattice space} in \S\ref{sec: transversal})
implies that the set $\mathcal{D}$  has zero Hausdorff dimension transversally to the orbit of this
$n$-dimensional group  (in particular, this set $\mathcal{D}$ has Hausdorff dimension~$n$; see Section
\ref{sec: definitions} and  Section \ref{sec: transversal}  for more details regarding  Hausdorff dimension
and tranversals, and 
\cite{Hee-TIFR}, \cite{Hee-unpublished} for more details regarding this application).

\Subsec{Measure rigidity}\label{intro measure}
The earliest results for  measure rigidity for
higher rank hyperbolic actions   deal with the Furstenberg problem:
\cite{RLyons}, \cite{Rudolph-2-and-3}, \cite{Johnson-invariant-measures}. Specifically,
Rudolph \cite{Rudolph-2-and-3} and Johnson \cite{Johnson-invariant-measures}
proved that  if
$\mu$ is a probability measure invariant and ergodic under the
action of the semigroup generated by $ \times m$, $ \times n$
(again with $m$, $n$ not powers of the same integer), and if some
element of this semigroup acts with positive entropy, then $\mu$
is Lebesgue.

When Rudolph's result appeared, the  second  author suggested
another test model for the measure  rigidity: two commuting
hyperbolic automorphisms of the three-dimensional torus. Since
Rudolph's proof  seemed, at least superficially,  too closely
related to symbolic dynamics, jointly with R.~Spatzier,  a more
geometric technique was developed. This allowed a unified
treatment of essentially all the classical examples of higher rank
actions for which rigidity of measures is expected
\cite{KatokSpatzier96}, \cite{KalininKatok99}, and in retrospect,
Rudolph's proof can also be interpreted in this framework.

This method (as well as most later work on measure rigidity for
these higher rank abelian actions) is based on the study of
conditional measures induced by a given invariant measure $ \mu$
on certain invariant foliations. The foliations considered include
stable and unstable foliations of various elements of the actions,
as well as intersections of such foliations, and are related to
the Lyapunov exponents of the action. For Weyl chamber flows these
foliations are given by orbits of unipotent subgroups normalized
by the action.

Unless there is an element of the action which acts with positive
entropy with respect to $\mu$, these conditional measures are
well-known to be $\delta$-measure supported on a single point, and
do not reveal any additional meaningful information about $\mu$.
Hence this and later techniques are limited to study actions where
at least one element has positive entropy. Under ideal situations,
such as the original motivating case of two commuting hyperbolic
automorphisms of the three torus, no further assumptions are
needed, and a result entirely analogous to Rudolph's theorem can
be proved using the method of \cite{KatokSpatzier96}.

However, for Weyl chamber flows, an additional assumption is
needed for the \cite{KatokSpatzier96} proof to work. This
assumption is satisfied, for example, if the flow along every
singular direction in the Weyl chamber is ergodic (though a weaker
hypothesis is sufficient). This additional assumption, which
unlike the entropy assumption is not stable under weak$^{*}$
limits, precludes applying the results  from
\cite{KatokSpatzier96} in many cases.

Recently, two new methods of proofs were developed, which overcome this difficulty.

The first method was developed by the first and second authors 
\cite{EinsiedlerKatok02}, following an idea mentioned at the end of 
\cite{KatokSpatzier96}. This idea uses the noncommutativity of the above-mentioned foliations (or more precisely, of the corresponding unipotent groups). This paper deals with general $\RR$-split semisimple Lie groups; in particular it is shown there that if
$\mu$ is an $A$-invariant measure on $X=\SL (k, \RR) / \Gamma$,
and if the entropies of $\mu$ with respect to {\it all} one-parameter groups are positive, then $\mu$ is the
Haar measure. It should be noted that for this method the properties of the lattice do not play any role, and
indeed this is true not only for $\Gamma = \SL (k, \ZZ)$ but for every discrete subgroup $\Gamma$. An
extension to the nonsplit case appeared in \cite{EinsiedlerKatokNonsplit}.
Using the methods we present in the second part of the present paper, the results of \cite{EinsiedlerKatok02} can be used to show that the set of exceptions to Littlewood's conjecture has Hausdorff dimension at most 1.

A different approach was developed by the third author, and was
used to prove a special case of the quantum unique ergodicity
conjecture \cite{Lindenstrauss03}. In its basic form, this
conjecture is related to the geodesic flow, which is not rigid, so
in order to be able to prove quantum unique ergodicity in certain
situations a more general setup for measure rigidity, following
Host \cite{Host-normal-numbers}, was needed. A special case of the
main theorem of \cite{Lindenstrauss03} is the following: Let  $A$
be an $\RR$-split Cartan subgroup of $\SL (2, \RR) \times \SL (2,
\RR)$.  Any $A$-ergodic measure on $\SL (2, \RR) \times \SL (2,
\RR) / \Gamma$ for which some one-parameter subgroup of $A$ acts
with positive entropy is algebraic. Here $\Gamma$ is e.g.\ an
irreducible lattice in $\SL (2, \RR) \times \SL (2, \RR)$. Since
the foliations under consideration in this case do commute, the
methods of \cite{EinsiedlerKatok02} are not applicable.

The method of \cite{Lindenstrauss03}
can be adapted to
quotients of more general groups, and in particular to $\SL (k, \RR)$.
It is noteworthy (and gratifying)
that for the space of lattices (and more general quotients of $\SL (k, \RR)$)
these two unrelated methods  are completely complementary: measures with
``high'' entropy (e.g.\ measures for which many one-parameter subgroup have positive entropy)
can be handled with the methods of \cite{EinsiedlerKatok02}, and measures with``low''
(but positive) entropy can be handled using the methods of 
\cite{Lindenstrauss03}.
Together, these methods give Theorem~\ref{theorem about lattice space}
(as well as the more general Theorem~\ref{t:slk} below for more general quotients).

The method of proof in \cite{Lindenstrauss03}, an adaptation of which we use
here, is based on study of the behavior of $\mu$ along certain unipotent
trajectories, using techniques introduced by Ratner in
\cite{Ratner-joinings}, \cite{Ratner-factors} to study unipotent flows, in particular
the H-property (these techniques are nicely exposed in Section~1.5 of   \cite{Witte-book}). This is
surprising because the techniques are applied on a measure $\mu$ which is {\it a~priori} not even
quasi-invariant under these (or any other) unipotent flows.

In showing that the high entropy and low entropy
cases are complementary we use a variant on the Ledrappier-Young entropy formula
\cite{LedrappierYoungII85}. Such use is one of the simplifying ideas in G.
Tomanov and Margulis' alternative proof of Ratner's theorem 
\cite{Margulis-Tomanov}.

\vskip6pt {\it Acknowledgment}. The authors are grateful to Dave Morris
Witte for pointing out some helpful references about nonisotropic
tori. E.L.\ would also like to thank Barak Weiss for introducing
him to this topic and for numerous conversations about both the 
Littlewood Conjecture and rigidity of multiparametric actions.
A.K.\ would like to thank Sanju Velani for helpful conversations
regarding the  Littlewood Conjecture. The authors would like to thank M. Ratner and the referees for many
helpful comments. The authors acknowledge the
hospitality of the Newton Institute for Mathematical Sciences in
Cambridge in the spring of 2000 and ETH Zurich in which some of the seeds of this
work have been sown. We would also like to acknowledge the hospitality of
the University of Washington, the Center for Dynamical Systems at  the Pennsylvania State University, and
Stanford University on more than one occasion.

\vglue10pt
\centerline{\bf Part I. Measure rigidity}
\vglue4pt 
Throughout this paper, let $G=\SL (k, \RR)$ for some $k \geq 3$, let $\Gamma$
be a discrete subgroup of $G$, and let $X = G / \Gamma$. As in the previous
section, we let $A<G$ denote the group of  $k \times k$ positive diagonal
matrices. We shall implicitly identify $$\Sigma=\{\bt \in \RR ^ k:t_1+\cdots
+t_k=0 \}$$ and the Lie algebra of $A$ via the map $(t _ 1, \dots, t
_ k) \mapsto \diag (t _ 1, \dots, t _ k)$. We write $\alpha ^\bt=\diag (e ^{t _
1}, \dots, e ^{t _ k})\in A$ and also $\alpha ^\bt$ for the left multiplication
by this element on $X$. This defines an $\RR ^{k-1}$ flow $\alpha$ on $X$.

A subgroup $U < G$ is {\it unipotent} if for every $g \in U$,  \
$g-\Id$ is nilpotent; i.e., for some $n$, \ $(g-\Id)^ n=0$. A group
$H$ is said to be {\it normalized by} $g \in G$ if $ g H g ^{-1} =
H$; $H$ is normalized by $L<G$ if it is normalized by every $g \in
L$; and the {\it normalizer} $N(H)$ of $H$ is the group of all $g
\in G$ normalizing it. Similarly, $g$ {\it centralizes} $H$ if
$gh=hg$ for every $h \in H$, and we set $C(H)$, the {\it
centralizer} of $H$ in $G$, to be the group of all $g \in G$
centralizing $H$.

If $U<G$ is normalized by $A$ then for every $x \in X$ and $a \in
A$, $a(Ux) = Uax$, so that the {\it foliation} of $X$ into $U$
orbits is invariant under the action of $A$. We will say that $a
\in A$ expands $U$ if all eigenvalues of $ \Ad (a)$ restricted to
the Lie algebra of $U$ are greater than one.

For any locally compact metric space $Y$ let $\mathcal{M} _ \infty (Y)$ denote
the space of Radon measures on $Y$ equipped with the weak$^*$ topology, i.e.\
all locally finite Borel measures on $Y$ with the coarsest topology for which
$\rho \mapsto \int _ Y f (y) d \rho (y)$ is continuous for every compactly
supported continuous $f$. For two Radon measures $\nu_1$ and $\nu_2$ on $Y$ we
write $$\nu_1 \propto \nu_2\mbox{ if }\nu_1=C \nu_2\mbox{ for some }C>0,
$$ and
say that $\nu_1$ and $\nu_2$ are proportional.

We let $B ^ Y_\varepsilon (y)$ (or $B _ \varepsilon (y)$ if $Y$ is understood) denote
the ball of radius $ \varepsilon$ around $y \in Y$; if $H$ is a group we set $B ^
H _ \varepsilon = B ^ H _ \varepsilon (\operatorname{I})$ where $\operatorname{I}$ is
identity in $H$; and if $H$ acts on $X$ and $x \in X$ we let $B ^ H _ \varepsilon
(x) = B ^ H _ \varepsilon \cdot x$.

Let $d(\cdot,\cdot)$ be the geodesic distance induced by a right-invariant
Riemannian metric on $G$. This metric on $G$ induces a right-invariant metric
on every closed subgroup $H \subset G$, and furthermore a metric on
$X=G/\Gamma$. These induced metrics we denote by the same letter.

\section{Conditional measures on $A$-invariant foliations,\\
 invariant measures, and
shearing}\label{sec: conditional}
\vglue-12pt 

\Subsec{Conditional measures}\label{sec: conditionalm}
A basic construction, which was introduced in the context of measure
rigidity in \cite{KatokSpatzier96} (and in a sense is already used implicitly
in \cite{Rudolph-2-and-3}), is the restriction of probability or even Radon
measures on a foliated space to the leaves of this foliation. A discussion can
be found in \cite[\S4]{KatokSpatzier96}, and a fairly general construction
is presented in \cite[\S3]{Lindenstrauss03}. Below we consider special cases of this general construction, summarizing its main
properties.

Let $ \mu$ be an $A$-invariant probability measure
on $X$. For any unipotent subgroup $U < G$ normalized by $A$, one
has a system  $ \left\{ \mu_{x,U} \right\} _ {x \in X}$ of Radon
measures on $U$ and a co-null set $X' \subset X$ with the
following properties\footnote{We are following the conventions of
\cite{Lindenstrauss03} in viewing the conditional measures
$\mu_{x,U}$ as measures on $U$. An alternative approach, which,
for example, is the one taken in \cite{KatokSpatzier96} and
\cite{KalininKatok99}, is to view the conditional measures as a
collection of measures on $X$ supported on single orbits of $U$;
in this approach, however, the conditional measure is not a Radon
measure on $X$, only on the single orbit of $U$ in the topology of
this submanifold.}:
\begin{enumerate}
\item The map $x \mapsto \mu _ {x, U}$ is measurable. \item For
every $ \varepsilon > 0$ and $x \in X '$, $\mu _ {x, U}
(B ^ U _ \varepsilon)>0$. \item  For every $x \in X'$ and $u \in U$
with $ux \in X '$, we have that $\mu _ {x,U} \propto (\mu _ {ux,
U}) u$, where  $(\mu _ {ux, U}) u$ denotes the push forward of the
measure $\mu _ {ux, U}$ under the map $v \mapsto vu$. \item For
every $\bt \in \Sigma$, and $x,\alpha ^\bt x \in X'$, $\mu _
{\alpha ^\bt x, U} \propto \alpha ^\bt (\mu _{x,U}) \alpha ^{-\bt}$.
\end{enumerate}
In general, there is no canonical way to normalize the measures $\mu _{x,U}$;
we fix a specific normalization by requiring that $\mu _{x,U} (B ^ U _ 1) = 1$
for every $x \in X'$. This implies the next crucial property.
\begin{enumerate}
\item[(5)] If $U \subset C(\alpha ^\bt)=\{ g \in G: g \alpha ^\bt=\alpha ^\bt g \}$
commutes with $\alpha ^\bt$, then
$\mu_ {\alpha ^\bt x, U} = \mu _{x,U}$
whenever $x,\alpha ^\bt x \in X'$.
\item[(6)] $\mu$ is $U$-invariant if, and only if, $\mu _ {x, U}$ is
a Haar measure on $U$ a.e.\ (see e.g.\ \cite{KatokSpatzier96} or the slightly more
general \cite[Prop.~4.3]{ Lindenstrauss03}).
\end{enumerate}

The other extreme to $U$-invariance occurs when $\mu _ {x, U}$ is atomic. If
$ \mu$ is $A$-invariant then outside some set of measure zero if $\mu _ {x, U}$
is atomic then it is supported on the identity $ \Id \in U$, in which case
we say that {\it $\mu _ {x, U}$ is trivial}. This follows from Poincar\'e
recurrence for an element $a \in A$ that uniformly expands the $U$-orbits (i.e.\
for which the $U$-orbits are contained in the unstable manifolds). Since the
set of $ x \in X$ for which $\mu _ {x, U}$ is trivial is $A$-invariant, if
$\mu$ is $A$-ergodic then either $\mu _ {x, U}$ is trivial a.s.\ or $\mu _ {x,
U}$ is nonatomic a.s. Fundamental to us is the following characterization of
positive entropy
(see \cite[\S~9]{Margulis-Tomanov} and \cite{KatokSpatzier96}):

\begin{enumerate}
\item[(7)] If for every $x \in X$ the orbit $Ux$ is the stable manifold through
$x$ with respect to $\alpha ^\bt$, then the measure theoretic entropy
$\h_\mu(\alpha ^\bt)$ is positive if and only if the conditional measures
$\mu_{x,U}$ are nonatomic a.e.
\end{enumerate}

So positive entropy implies that the conditional measures are nontrivial a.e.,
and the goal is to show that this implies that they are Haar measures. Quite
often one shows first that the conditional measures are translation invariant
under some element up to proportionality, which makes the following observation useful.
\begin{enumerate}
\item[(8)] Possibly after replacing $X'$ of (1)--(4) by a conull subset,
we see that for any
$x \in X'$ and any $u \in U$ with $\mu_{x,U}\propto \mu_{x,U}u$, in fact,
$\mu_{x,U}=\mu_{x,U}u$ holds.
\end{enumerate}
This was first   shown in \cite{KatokSpatzier96}.
The proof of this fact only uses Poincar\'e
recurrence and (4) above; for completeness we provide a proof below.

\demo{Proof of {\rm (8)}} Let $\bt$ be such that $\alpha ^\bt$ uniformly contracts the
$U$-leaves (i.e.\ for every $x$ the $U$-orbit $Ux$ is part of the stable
manifold with respect to $\alpha ^\bt$). Define for $M>0$
$$
D_M = \left\{ x \in X': \mu_{x,U}\bigl(B_2 ^ U\bigr)<M \right\}
.$$ We claim that for every $x \in X' \cap \bigcup_ M \limsup_{n
\to \infty} \alpha ^{-n \bt} D_M$ (i.e. any $x \in X'$ so that $ \alpha ^ {n
\bt}$ is in $D_M$ for some $M$ for infinitely many $n$) if $\mu_{x,U} = c
\mu_{x,U}u$ then $c \leq 1$.

Indeed, suppose $x \in X' \cap \limsup _{n \to \infty} \alpha ^{-n \bt} D_M$ and $u \in U$ satisfy $\mu_{x,U} = c \mu_{x,U}u$. Then for any $n,k$
$$
\mu_{\alpha ^ {n \mathbf t} x,U} = c^k \mu_{\alpha ^ {n \mathbf t} x,U} (\alpha ^ {n \mathbf t} u^k \alpha ^ {- n \mathbf t})
.$$
Choose $k>1$ arbitrary. Suppose $n$ is such that $\alpha ^ {n \mathbf t} x \in D_M$ and suppose that $n$ is sufficiently large that $\alpha ^ {n \mathbf t} u^k \alpha ^ {- n \mathbf t} \in B_1^U$, which is possible since $\alpha ^ {\mathbf t}$ uniformly contracts $U$. Then
\begin{align*}
M & \geq \mu _ {\alpha ^ {n \mathbf t}x,U} (B_2^U) \geq
\mu_{\alpha ^ {n \mathbf t}x,U} (B_1^U \alpha ^ {n \mathbf t} u^{k} \alpha ^ {- n \mathbf t})\\
& = (\mu_{\alpha ^ {n \mathbf t} x,U} \alpha ^ {n \mathbf t} u^{-k} \alpha ^ {- n \mathbf t})(B_1^U) \\
& =
c^k \mu_{\alpha ^ {n \mathbf t} x,U} (B_1^U) = c^k
.\end{align*}
Since $k$ is arbitrary this implies $c \leq 1$.

If $\mu_{x,U} = c \mu_{x,U}u$ then $\mu_{x,U} = c^{-1} \mu_{x,U}u^{-1}$, so the above argument applied to $u ^{-1}$ shows that $c \geq 1$, hence $\mu_{x,U} = \mu_{x,U}u$.

Thus we see that if we replace $X'$ by  $X' \cap \bigcup_M\limsup _{n \to
\infty} \alpha ^{-n \bt} D_M$ --- a conull subset of $X'$,  then (8) holds
for any $x \in X'$.
\Endproof\vskip4pt  

Of particular importance to us will be the following one-parameter
unipotent subgroups of $G$, which are parametrized by pairs $(i,j)$
of distinct integers in the range $\left\{ 1, \dots,k \right\}$:
$$
u_{ij} (s)=\exp(s E_{ij})=\Id+s E_{ij}, \qquad U_{ij}=\{ u _
{ij}(s):s \in \RR \},
$$
where $E_{ij}$ denotes the matrix with 1 at the $i^{\rm th}$ row and
$j^{\rm th}$ column and zero everywhere else. It is easy to see that  these
groups are normalized by $ A$; indeed, for $ \mathbf t = (t _ 1,
\dots, t _ k) \in \Sigma$
$$
\alpha ^ \mathbf t  u _ {ij}(s) \alpha ^ {-\mathbf t} =   u _
{ij}(e ^{t_i - t_j}s) .
$$
Since these groups are
normalized by $A$, the orbits of $U _ {ij}$ form an $A$-invariant foliation of
$X = \SL (k, \RR) / \Gamma$ with one-dimensional leaves. We will use $\mu ^
{ij} _ x$ as a shorthand for $\mu _ {x, U _ {ij}}$; any integer $i \in \left\{
1, \dots,k \right\}$ will be called an index; and unless otherwise stated, any
pair $i,j$ of indices is implicitly assumed to be distinct.

Note that for the conditional measures $\mu_x ^{ij}$ it is easy to find a
nonzero $\bt \in \Sigma$ such that (5) above holds; for this all we need is
$t_i=t_j$. Another helpful feature is the one-dimensionality of $U_{ij}$ which
also helps to show that $\mu_x ^{ij}$ are a.e.\ Haar measures.
In particular we
have the following:
\begin{enumerate}
\item[(9)] Suppose there exists a set of positive measure $B \subset X$ such
that for any $x \in B$ there exists a nonzero $u \in U_{ij}$ with
$\mu_x ^{ij}\propto \mu_{x}^{ij}u$. Then for a.e.\ $x \in B$ in fact $\mu_x ^{ij}$
is a Haar measure of $U_{ij}$, and if $\alpha$ is ergodic then $\mu$ is
invariant under $U_{ij}$.
\end{enumerate}

{\it Proof of {\rm (9)}}.
Recall first that by (8) we can assume
$\mu_x ^{ij}=\mu_{x}^{ij}u$ for $x \in B$. Let $K \subset B$ be a compact set of
measure almost equal to $\mu(B)$ such that $\mu_x ^{ij}$ is continuous for $x \in
K$. It is possible to find such a $K$ by Luzin's theorem. Note however, that
here the target space is the space of Radon measures $\mathcal
M_\infty(U_{ij})$ equipped with the weak$^*$ topology so that a more general
version \cite[p.\ 69]{Fed69} of Luzin's theorem is needed. Let $\bt \in \Sigma$
be such that $U_{ij}$ is uniformly contracted by $\alpha ^\bt$. Suppose now
$x \in K$ satisfies Poincar\'e recurrence for every neighborhood of $x$ relative
to $K$. Then there is a sequence $x_\ell=\alpha ^{n_\ell\bt}\in K$ that
approaches $x$ with $n_\ell \rightarrow \infty$. Invariance of $\mu_x ^{ij}$ under
$u$ implies invariance of $\mu_{x_\ell}$ under the much smaller element
$\alpha ^{n_\ell\bt} u \alpha ^{-n_\ell\bt}$ and all its powers. However, since
$\mu_{x_\ell}^{ij}$ converges to $\mu_x ^{ij}$ we conclude that $\mu_x ^{ij}$ is
a Haar measure of $U_{ij}$. The final statement follows from (4) which implies
that the set of $x$ where $\mu_x ^{ij}$ is a Haar measure is $\alpha$-invariant.
\Endproof\vskip4pt  

Even when $\mu$ is not invariant under $U_{ij}$ we still have the following
maximal ergodic theorem \cite[Thm.\ A.1]{Lindenstrauss03} proved by the last
named author in joint work with  D.~Rudolph, which is related to a maximal ergodic theorem of
Hurewicz~\cite{Hurewicz}.
\begin{enumerate}
\item[(10)] For any $f \in L ^ 1(X,\mu)$ and $\alpha >0$,  
$$
\mu\Bigl(\Bigl\{ x: \int_{B_r ^{U_{ij}}} f(ux)\operatorname{d}\!\mu_x ^{ij}> \alpha
\mu_x ^{ij}\bigl(B_r ^{U_{ij}}\bigr)\mbox{ for some
}r>0\Bigr\}\Bigr)<\frac{C \| f \|_1}{\alpha}
$$
for some universal constant $C>0$.
\end{enumerate}

\vglue-8pt
\Subsec{Invariant measures{\rm ,} high and low entropy cases}
We are now in a position to state the general measure rigidity result for quotients of $G$:

\begin{theorem}\label{t:slk}
Let $X= G / \Gamma$ and $A$ be as above. Let $\mu$ be an
$A$-invariant and ergodic probability measure on $X$. For any pair
of indices $a,b${\rm ,} one of the following three properties must hold.
\begin{enumerate}
\item  The conditional measures $\mu _ {x}^{ab}$ and $\mu_{x}^{ba}$ are trivial
a.e. \item  The conditional measures $\mu _ {x}^{ab}$ and $\mu_{x}^{ba}$ are
Haar a.e.{\rm ,} and $\mu$ is invariant under left multiplication with elements of
$H_{ab}=\langle U_{ab},U_{ba}\rangle$. \item  \label{exceptional returns} Let
$A _ {ab}' = \left\{ \alpha ^ \mathbf s: \mathbf s \in \Sigma \text{ and }
s_a=s_b \right\}$. Then a.e.\ ergodic component of $\mu$ with respect to $A' _
{ab} $ is supported on a single $C(H_{ab})$-orbit{\rm ,} where $C(H_{ab})=\{ g \in G:
gh=hg$ for all $h \in H_{ab}\}$ is the centralizer of $H_{ab}$.
\end{enumerate}
\end{theorem}

{\it Remark}. If $k=3$
then (3) is equivalent to the following:
\begin{enumerate}
\item[(3$'$)] There exist  a nontrivial $\bs \in \Sigma$ with $s_a=s_b$  and a
point $x_0 \in X$ with $\alpha ^\bs x_0=x_0$ such that the measure $\mu$ is
supported by the orbit of $x_0$ under $C(A_{ab}')$. In particular, a.e.\ point
$x$ satisfies $\alpha ^\bs x=x$.
\end{enumerate}
Indeed, in this case
$C(H_{ab})$ contains only diagonal matrices, and Poincar\'e recurrence for
$A_{ab}'$ together with (3) imply that a.e.\ point is periodic under $A_{ab}'$.
However, ergodicity of $\mu$ under $A$ implies that the period $\bs$ must be
the same a.e. Let $x_0 \in X$ be such that every neighborhood of $x_0$ has
positive measure. Then $x$ close to $x_0$ is fixed under $\alpha ^\bs$ only if
$x \in C(A_{ab}')x_0$, and ergodicity shows (3$'$).
The examples of M.~Rees \cite{Rees}, \cite[\S9]{EinsiedlerKatok02} of nonalgebraic $A$-ergodic measures in certain quotients of $\SL (3, \RR)$ (which certainly can have positive entropy) are precisely of this form, and show that case (3) and 
(3$'$)
above are not superfluous.

When $ \Gamma = \SL (k, \ZZ)$, however, this
phenomenon, which we term {\it exceptional returns}, does not happen. We will
show this in Section \ref{s: exceptional returns}; similar observations
have been made earlier in \cite{Margulis-conjectures},
\cite{Lindenstrauss-Weiss}. We also refer the reader to \cite{Tomanov-maxtori} for a
treatment of similar questions for inner lattices in $ \SL (k, \RR)$ (a certain class
of 
lattices in $\SL (k, \RR)$).

The conditional measures $ \mu _ {x}^{ij}$ are intimately
connected with the entropy. More precisely, $\mu$ has positive
entropy with respect to $ \alpha ^ \mathbf t$ if and only if for
some $i,j$ with $ t_i > t_j$ the measures $\mu ^{ij} _ {x}$ are not
a.s.\ trivial (see Proposition~\ref{p: ent_cond} below for more
details; this fact was first proved in \cite{KatokSpatzier96}).
Thus (1) in Theorem~\ref{t:slk} above holds for all pairs of
indices $i,j$ if, and only if, the entropy of $\mu$ with respect
to every one-parameter subgroup of $A$ is zero.

In order to prove Theorem~\ref{t:slk}, it is enough to show that
for every $a,b$ for which the $\mu _ {x}^{ab}$ is a.s.\ nontrivial
either Theorem~\ref{t:slk}.(2) or Theorem~\ref{t:slk}.(3) holds.
For each pair of indices $a,b$, our proof is divided into two
cases which we loosely refer to as the high entropy and the low
entropy case:

\demo{High entropy case} There is an additional pair of indices $i,j$
distinct from $a,b$ such that $i=a$ or $j=b$  for which $\mu _ {x}^{ij}$ are
nontrivial a.s. In this case we prove:

\begin{theorem} \label{t: high entropy}
If both $\mu _ {x}^{ab}$ and $\mu _ x ^ {ij}$ are nontrivial a.s.{\rm ,} for distinct
pairs of indices $i,j$ and $a,b$ with either $i=a$ or $j=b${\rm ,} then both $\mu _
x ^{ab}$ and $\mu _ x ^ {ba}$ are in fact Haar measures a.s.\ and $\mu$ is
invariant under $H_{ab}$.
\end{theorem}

The proof in this case, presented in Section~\ref{s: high entropy} makes
use of the noncommutative structure of certain unipotent
subgroups of $G$, and follows  \cite{EinsiedlerKatok02} closely.
However, by careful use of an adaptation of a formula of
Ledrappier and Young (Proposition~\ref{p: ent_cond} below)
relating entropy to the conditional measures $\mu _ {x}^{ab}$ we
are able to extract some additional information. It is interesting
to note that Margulis and Tomanov used the Ledrappier-Young theory
for a similar purpose in \cite{Margulis-Tomanov}, simplifying
some of Ratner's original arguments in the classification of
measures invariant under the action of unipotent groups.

\demo{Low entropy case} For every pair of indices $i,j$ distinct from $a,b$
such that $i=a$ or $j=b$,  $\mu _ {x}^{ij}$ are trivial a.s. In this case there
are two possibilities:

\begin{theorem} \label{t: low entropy}
Assume $\mu _ x ^{ab}$ are a.e.\ nontrivial{\rm ,} and $\mu _ {x}^{ij}$ are trivial
a.e.\ for every pair $i,j$ distinct from $a,b$ such that $i=a$ or $j=b$. Then
one of the following properties holds.
\begin{enumerate}
\item $\mu$ is $U _ {ab}$ invariant.
\item Almost every $A'_{ab}$-ergodic component of $\mu$ is
supported on a single $C(H _ {ab})$ orbit.
\end{enumerate}
\end{theorem}

We will see in Corollary \ref{cor: low symmetry}
that in the low entropy case $\mu _ {x} ^ {ba}$ is also
nontrivial; so applying Theorem~\ref{t: low entropy} for $ U _
{ba}$ instead of $U _ {ab}$ one sees that either $\mu$ is $ H _
{ab}$-invariant or almost every $A'_{ab}$-ergodic component of
$\mu$ is supported on a single $C(H _ {ab})=C(H_{ba})$ orbit.

In this case we employ the techniques developed by the third named author in
\cite{Lindenstrauss03}. There, one considers invariant measures on irreducible
quotients of products of the type $ \SL (2, \RR) \times L$ for some algebraic
group $L$. Essentially, one tries to prove a Ratner type result (using methods
quite similar to Ratner's \cite{Ratner-factors}, \cite{Ratner-joinings}) for the $U _
{ab}$ flow even though $\mu$ is not assumed to be invariant or even quasi
invariant under $U _ {ab}$. Implicitly in the proof we use a variant of
Ratner's H-property (related, but distinct from the one used by Witte in
\cite[\S6]{Witte-rigidity}) together with the maximal ergodic theorem  for
$U_{ab}$ as in (9) in Section~\ref{sec: conditionalm}.

\section{More about entropy and the high entropy case}\label{s: high entropy}

A well-known theorem by Ledrappier and Young
\cite{LedrappierYoungII85} relates the entropy, the dimension of
conditional measures along invariant foliations, and Lyapunov
exponents, for a general $C ^ 2$ map on a compact manifold, and in
\cite[\S9]{Margulis-Tomanov} an adaptation of the general
results to flows on locally homogeneous spaces is provided. In the
general context, the formula giving the entropy in terms of the
dimensions of conditional measures along invariant foliations
requires consideration of a sequence of subfoliations, starting
from the foliation of the manifold into stable leaves. However,
because the measure $ \mu$ is invariant under the full $A$-action
one can relate the entropy to the conditional measures on the
one-dimensional foliations into orbits of $U_{ij}$ for all pairs
of indices $i,j$.

We quote the following from \cite{EinsiedlerKatok02}; in that paper, this
proposition is deduced from the fine structure of the conditional measures on
full stable leaves for\break $A$-invariant measure; however, it can also be deduced
from a more general result of Hu regarding properties of commuting
diffeomorphisms \cite{Hu94}. It should be noted that the constants $s_{ij}
(\mu)$ that appear below have explicit interpretation in terms of the pointwise
dimension of $ \mu _ {x}^{ij}$ \cite{LedrappierYoungII85}.

\begin{proposition}[{\cite[Lemma 6.2]{EinsiedlerKatok02}}]\label{p: ent_cond}
Let $\mu$ be an $A$-invariant and ergodic probability measure on $X = G / \Gamma$ with $G = \SL (k, \RR)$ and $\Gamma < G$ discrete.
Then for any pair of indices $i,j$ there are constants $s_{ij}(\mu) \in [0,1]$ so that\/{\rm :}\/
\begin{enumerate}
\item $s_{ij} (\mu) = 0$ if and only if for a.e.\ $x${\rm ,} $ \mu _ x ^{ij}$ are
atomic and supported on a single point. 
\item If a.s.\ $\mu _ x ^{ij}$ are Haar
\/{\rm (}\/i.e.\ $\mu$ is $U _ {ij}$ invariant\/{\rm ),}\/ then $s_{ij} (\mu) =1.$ 
\item For any
$\bt \in \Sigma,$
\begin{equation}\label{e:LedYoung}
\h_ \mu(\alpha ^\bt)=\sum_{i,j} s_{ij} (\mu) (t_i - t_j)^+.
\end{equation}
\end{enumerate}
\end{proposition}

Here $(r)^+=\max(0,r)$ denotes the positive part of $r \in \RR$.

We note that the converse to (2) is also true. A similar proposition holds for more general semisimple groups $G$.
In particular we get the following (which is also proved in a somewhat different way in \cite{KatokSpatzier96}):

\begin{corollary}
For any $\mathbf t \in \Sigma${\rm ,} the entropy $\h_ \mu(\alpha ^\bt)$ is positive
if and only if there is a pair of indices $i,j$ with $t_i - t_j>0$ for which $
\mu _ x ^{ij}$ are nontrivial a.s.
\end{corollary}

A basic property of the entropy is that for any $\mathbf t \in \Sigma$,
\begin{equation} \label{e: past and future}
\h _ \mu (\alpha ^ \mathbf t) = \h _ \mu (\alpha ^ {- \mathbf t})
.\end{equation} As we will see this gives nontrivial identities between the $ s
_ {ij} (\mu)$. 

\begin{figure}[ht]
\begin{center}
 \epsfig{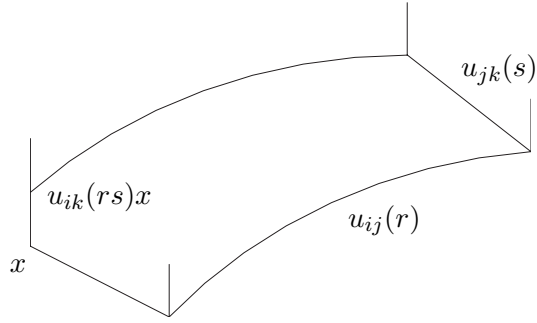}
\setlength{\unitlength}{1mm}
\begin{picture}(0,0)

\ifpdf \put(-90,0){ \put(15,9){$x$} \put(20,18){$u_{ik}(rs)x$}
\put(60,15){$u_{ij}(r)$} \put(77,38){$u_{jk}(s)$} }
\else \put(-87,-3){ \put(15,9){$x$}
\put(20,18){$u_{ik}(rs)x$} \put(60,15){$u_{ij}(r)$} \put(75,35){$u_{jk}(s)$} }
\fi
\end{picture}
\end{center}
\caption{\label{EKfigure} One key ingredient of the proof of
Lemma \ref{l:commEK} in \cite{EinsiedlerKatok02} is the translation produced along $U_{ik}$
when going along $U_{ij}$ and $U_{jk}$ and returning to the same
leaf $U_{ik}x$.}
\end{figure}

The following is a key lemma from \cite{EinsiedlerKatok02}; see Figure
\ref{EKfigure}.

\begin{lemma}[{\cite[Lemma 6.1]{EinsiedlerKatok02}}]\label{l:commEK}
Suppose $\mu$ is an $A$-invariant and ergodic probability measure{\rm ,}
 $i,j,k$ distinct indices  such that both
$\mu_x ^{ij}$ and $\mu_x ^{jk}$ are non\-atomic a.e.
Then $\mu$ is $U _ {ik}$-invariant.
\end{lemma}

\demo{Proof of Theorem~{\rm \ref{t: high entropy}}}
For $\ell=a,b$ we define the sets
\begin{eqnarray*}
C_\ell&=&\{ i \in \{ 1,\dots,k \} \setminus \left\{ a,b \right\}:s_{i \ell}(\mu)>0 \},\\
R_\ell&=&\{ j \in \{ 1,\dots,k \} \setminus \left\{ a,b \right\}:s_{\ell j}(\mu)>0 \},\\
C_\ell ^ L&=&\{ i \in \{ 1,\dots,k \} \setminus \left\{ a,b \right\}: \text{{$\mu$} is $U_{i
\ell}$-invariant}\},\\ R_\ell ^ L&=&\{ j \in \{ 1,\dots,k \} \setminus \left\{ a,b \right\}:
\text{{$\mu$} is $U_{\ell j}$-invariant}\} .
\end{eqnarray*}
Suppose $i \in C_a$; then the conditional measures $\mu_x ^{ia}$
are nontrivial a.e.\ by Proposition~\ref{p: ent_cond}. Since by
assumption $\mu_x ^{ab}$ are nontrivial a.e., Lemma \ref{l:commEK}
shows that $\mu_x ^{ib}$ are Lebesgue a.e. This shows that $C_a
\subset C_b ^ L$, and $R_b \subset R_a ^ L$ follows similarly.

Let $\bt=(t_1, \dots, t _ k)$ with $t _ i = -1/k$ for $i \ne a$ and $t _ a =
1-1/k$. For the following expression set $s_{aa}=0$. 
By Proposition~\ref{p: ent_cond} the entropy of $\alpha
^{\bt}$ equals
\begin{eqnarray}\label{e:ent_alphas}
\h_ \mu(\alpha ^{\bt})&=&s_{a1} (\mu) +\cdots +s_{ak} (\mu)\\
&=&s_{ab}(\mu)+|R_a ^ L|+\sum_{j \in R_a \setminus R_a ^ L}s_{aj}(\mu)>
|R_a ^ L|,\nonumber
\end{eqnarray}
where we used our assumption that $s_{ab}(\mu)>0$. Applying
Proposition~\ref{p: ent_cond} for $\alpha ^{-\bt}$ we see similarly
that
\begin{equation}\label{e:ent_alphas2}
\h_ \mu(\alpha ^{-\bt})=s_{1a}(\mu)+\cdots
+s_{ka}(\mu)=s_{ba}(\mu)+\sum_{i \in C_a} s_{ia}(\mu) \leq
(1+|C_a|),
\end{equation}
where we used the fact   that $s_{ia}(\mu) \in[0,1]$ for $a=2,\dots,k$.
However, since the entropies of $\alpha ^\bt$ and of
$\alpha ^{-\bt}$ are equal, we get $|R_a ^ L|\leq |C_a|$.

Using $\bt'=(t _ 1 ', \dots, t _ k ')$ with $t _ i' = - 1/k$ for
$i \ne b$ and $t _ b' = 1-1/k$ instead of $\bt$ in the above
paragraph shows similarly $|C_b ^ L|\leq |R_b|$. Recall that $C_a
\subset C_b ^ L$ and $R_b \subset R_a ^ L$. Combining these
inequalities we conclude that
$$
|R_a ^ L|\leq |C_a|\leq |C_b ^ L|\leq |R_b|\leq |R_a ^ L|,
$$
and so all of these sets have the same cardinality. However, from
\eqref{e:ent_alphas} and \eqref{e:ent_alphas2} we see that $s_{ab}(\mu)+|R_a ^
L|\leq \h_\mu(\alpha ^\bt)\leq s_{ba}(\mu)+|C_a|$. Together we see that
\begin{equation}\label{e:sbasab}
 s_{ba}(\mu)\geq s_{ab}(\mu)>0.
\end{equation} 
From this we conclude as before that $C_a
\subset C_b ^ L \subset C_a ^ L$, and so $C_a=C_a ^ L$. Similarly, one sees
$R_b=R_b ^ L$.

This shows that if $s_{ab} (\mu) > 0$ and $s _ {ij} (\mu)>0$ for some other
pair $i,j$ with either $i=a$ or $j = b$, then in fact $\mu$ is $U _
{ij}$-invariant. If there was at least one such pair of indices $i,j$ we could
apply the previous argument to $i,j$ instead of $a,b$ and get that $\mu$ is
$U_{ab}$-invariant. 
\Endproof\vskip4pt  

In particular, we have seen in the proof of Theorem \ref{t: high
entropy} that $s_{ab}>0$ implies \eqref{e:sbasab}.  We conclude the following symmetry.

\begin{corollary}\label{cor: low symmetry}
For any pair of indices $(a,b)${\rm ,}  $s_{ab}=s_{ba}$ . In particular{\rm ,} $\mu _
{x}^{ab}$ are nontrivial a.s.{\rm ,} if and only if{\rm ,} $\mu _ {x}^{ba}$ are nontrivial
a.s.
\end{corollary}

\section{The low entropy case}\label{s: low entropy}

We let $A'_{ab}=\{ \alpha ^\bs \in A:s_a=s_b \}$, and let $\alpha ^{\mathbf s}
\in A ' _ {ab}$. Then $\alpha ^ {\mathbf s}$ commutes with $U _ {ab}$, which
implies that $\mu_x ^{ab}=\mu_{\alpha ^\bs x}^{ab}$ a.e.

For a given pair of indices $a, b$, we define the following subgroups of $G$:
\begin{eqnarray*}
L _ {(ab)}& = & C(U_{ab}), \\
U _ {(ab)}& =& \langle U _ {ij}: \text{{$i = a$} or $j = b$} \rangle,\\
C _ {(ab)} &=& C (H _ {ab}) = C(U_{ab})\cap C(U_{ba}).
\end{eqnarray*}

Recall that the metric on $X$ is induced by a right-invariant
metric on $G$. So for every two $x,y \in X$ there exists a $g \in G$
with $y=gx$ and $d(x,y)=d(\Id,g)$.

\Subsec{Exceptional returns}\label{sec: exceptional}

\begin{definition}
We say for $K \subset X$ that the {\em $A' _ {ab}$-returns to $K$
are exceptional \/{\rm (}\/strong exceptional\/{\rm )}\/} if there exists a $\delta>0$ so that for
all $x, x' \in K$, and $\alpha ^\bs \in A'_{ab}$ with
$x'=\alpha ^{\bs}x \in B_ \delta(x)\cap K$  every $g
\in B_\delta ^ G$ with $x'=gx$ satisfies $g \in L_{(ab)}$ ($g \in
C_{(ab)}$ respectively).
\end{definition}

\begin{lemma}\label{l:exceptional}
There exists a null set $N \subset X$ such that for any compact
$K \subset X \setminus N$ with exceptional $A'_{ab}$-returns to $K$
the $A'_{ab}$-returns to $K$ are in fact strong exceptional.
\end{lemma}

\Proof 
To simplify notation, we may assume without loss of generality that $a=1,b=2$, and write $A'$, $U$, $L$,
$C$ for
$A'_{12}$, $U_{(12)}$, $L _ {(1 2)}$, $C _ {(1 2)}$ respectively. We write, for a given matrix $g \in G$,
\begin{equation}\label{e:matrixA}
g= \left(
{\begin{array}{ccc} a_1 & g_{12} & g_{1*}\\
g_{21} & a_2 & g_{2*}  \\
g_{*1} & g_{*2} & a_*
\end{array}}
\right),
\end{equation}
with the understanding that $a_1,a_2,g_{12},g_{21} \in \RR $,
$g_{1*}$, $g_{2*}$ (resp.\ $g_{*1}$, $g_{*2}$) are row (resp.\
column) vectors with $k-2$ components, and $a_*\in \Mat(k-2,\RR)$.
(For $k=3$ of course all of the above are real numbers, and we can
write $3$ instead of the symbol $*$.) Then $g \in L$ if and only
if $a_1 = a_2$ and $g_{21}$, $g _ {{*} 1}$, $g _ {2 {*}}$ are all
zero.  $g \in C$ if in addition $g _ {1 2}$,  $g _ {1 {*}}$, $g _
{{*} 2}$ are zero.

For $\ell \geq 1$ let $D_\ell$ be the set of $x \in X$ with the
property that for all $z \in B_{1/\ell}(x)$ there exists a unique
$g \in B ^ G_{1/\ell}$ with $z=gx$. Note that
$\bigcup_{\ell=1}^ \infty D_\ell=X$, and that for every compact set,
$K \subset D_\ell$ for some $\ell>0$.

Let first $\alpha ^\bs \in A'$ be a fixed element, and let
$E_{\ell,\bs}\subset D_\ell$ be the set of points $x$ for which
$x'=\alpha ^{\bs}x \in B_{1/\ell}(x)$ and $x' = gx$ with $g \in
B ^ G_{1/\ell} \cap L = B ^ L_{1/\ell}$. Since $g \in B ^ G_{1/\ell}$ is
uniquely determined by $x$ (for a fixed $\bs$), we can define (in
the notation of \eqref{e:matrixA}) the measurable function
$$
f(x)=\max\bigl(|g_{12}|,\| g_{1*}\|,\| g_{*2}\|\bigr)\mbox{ for
}x \in E_{\ell,\bs}.
$$

Let $\bt=(-1,1,0,\dots,0)\in \Sigma$. Then conjugation with
$\alpha ^\bt$ contracts $U$. In fact for $g$ as in
\eqref{e:matrixA} the entries of $\alpha ^{\bt}g \alpha ^{-\bt}$
corresponding to $g_{12},g_{1*}$ and $g_{2*}$ are
$e ^{-2}g_{12},e ^{-1}g_{1*}$ and $e ^{-1}g_{2*}$, and those
corresponding to $g_{21},g_{*1}$ and $g_{*2}$ are
$e ^{2}g_{21},eg_{*1}$ and $eg_{*2}$. Notice that the latter are
assumed to be zero. This shows that for $x \in E_{\ell,\bs}$ and
$\alpha ^{-n\bt}x \in D_\ell$, in fact 
$\alpha ^{-n\bt}x \in E_{\ell,\bs}$. Furthermore $f(\alpha ^{-n\bt}x)\leq
e ^{-n}f(x)$. Poincar\'e recurrence shows that $f(x)=0$ for a.e.\ $x \in
E_{\ell,\bs}$ -- or equivalently $\alpha ^\bs x \in B_{1/\ell}^ C(x)$ for a.e.\
$x \in D_\ell$ with $\alpha ^\bs x \in B ^ L_{1/\ell}(x)$.

Varying $\bs$ over all elements of $\Sigma$ with rational
coordinates and $\alpha ^\bs \in A'$, we arrive at a nullset
$N_\ell \subset D_\ell$ so that $\alpha ^{\bs}x \in B ^ L_{1/\ell}(x)$
implies $\alpha ^ {\bs} x \in B ^ C_{1/\ell}(x)$ for all such
rational $\bs$. Let $N$ be the union of $N_{\ell}$ for
$\ell=1,2,\dots\  $. We claim that $N$ satisfies the lemma.

So suppose $K \subset X \setminus N$ has $A'$-exceptional returns. Choose $\ell
\geq 1$ so that $K \subset D_{\ell}$, and furthermore so that $\delta=1/\ell$
can be used in the definition of\break $A'$-exceptional returns to $K$. Let $x \in
K$, $x'=\alpha ^\bs x \in B_{1/\ell}(x)$ for some $\bs \in \Sigma$ with
$\alpha ^\bs \in A'$, and $g \in B_{1/\ell}^ G$ with $x'=gx$. By assumption on
$K$, we have that $g \in L $. Choose a rational $\tilde\bs \in \Sigma$ close to
$\bs$ with $\alpha ^{\tilde\bs}\in A'$ so that $\alpha ^{\tilde\bs}x \in
B_{1/\ell}(x)$. Clearly $\tilde g=\alpha ^{\tilde\bs-\bs}g$ satisfies
$\alpha ^{\tilde\bs}x=\tilde gx$ and so $\tilde g \in B_{1/\ell}^ L$. Since $x
\in K \subset D_{1/\ell}\setminus N_{1/\ell}$, it follows that $\tilde g \in
C$. Going back to $x'=\alpha ^\bs x$ and $g$ it follows that $g \in C$.
\Endproof\vskip4pt  

Our interest in exceptional returns is explained by the following
proposition. Note that condition (1) below is exactly Theorem
\ref{t: low entropy}(2).

\begin{proposition}\label{p: equivalence}
For any pair of indices $a,b$ the following two conditions are
equivalent. \vglue-20pt \phantom{up}
\begin{enumerate}
\item A.e.\ ergodic component of $\mu$ with respect to $A_{ab}'$
is supported on a single $C_{(ab)}$-orbit. \item For every
$\varepsilon>0$ there exists a compact set $K$ with measure
$\mu(K)>1-\varepsilon$ so that the $A'_{ab}$-returns to $K$ are strong
exceptional.
\end{enumerate}
\end{proposition}
\phantom{up}
\vglue-19pt

The ergodic decomposition of $\mu$ with respect to $A '_{ab}$ can
be constructed in the following manner: Let $\mathcal E'$ denote
the $\sigma$-algebra of Borel sets which are $A'_{ab}$ invariant.
For technical purposes, we use the fact that $(X, \cB_X, \mu)$ is
a Lebesgue space to replace $\mathcal E'$ by an equivalent
countably generated sub-sigma algebra $\mathcal{E}$. Let
$\mu_x ^{\mathcal E}$ be the family of conditional measures of
$\mu$ with respect to the $\sigma$-algebra $\mathcal E$. Since
$\mathcal{E}$ is countably generated the atom $[x] _
{\mathcal{E}}$ is well defined for all $x$, and it can be arranged
that for {\it all} $x$ and $y$ with $y \in [x] _ {\mathcal{E}}$
the conditional measures $\mu _ x ^ \mathcal{E} = \mu _ y ^
\mathcal{E}$, and that for all $x$, $\mu _ x ^ \mathcal{E}$ is a
probability measure.

Since $\mathcal E$ consists of $A '_{ab}$-invariant sets, a.e.\
conditional measure is $A '_{ab}$-invariant, and can be shown to
be ergodic. So the decomposition of $\mu$ into conditionals
\begin{equation}\label{e: decomposition}
\mu=\int_X \mu_x ^{\mathcal E}\operatorname{d}\!\mu
\end{equation}
gives the ergodic decomposition of $\mu$ with respect to $A
'_{ab}$.

\Proof 
For simplicity, we write $A ' = A ' _ {ab}$ and $C = C _ {(ab)}$.

\demo{$(1) \Longrightarrow (2)$}
Suppose a.e.\ $A '$ ergodic component is supported on a single
$C$-orbit. Let $\varepsilon>0$. For any fixed $r>0$ we define
$$
f_r(x)=\mu_x ^{\mathcal E}(B_r ^ C(x)).
$$
By the assumption $f_r(x)\nearrow 1$ for $r \rightarrow \infty$
and a.e.\ $x$. Therefore, there exists a fixed $r>0$ with
$\mu(C_r)>1-\varepsilon$, where $C_r=\{ x:f_r(x)>1/2 \}$.

Fix some $x \in X$. We claim that for every small enough $\delta>0$
\begin{equation}\label{e: intersection equation}
B ^ C _ {2 r} (x) \cap B _ \delta (x) =B ^ C _ \delta (x).
\end{equation}
Indeed, by the choice of the metric on $X$ there exists
$\delta'>0$ so that the map $g \mapsto gx$ from $B_{3 \delta'}^ G$
to $X$ is an isometry. Every $g \in B_{2r}^ C$ satisfies that
either $d(B_{\delta'}^ C(x),B_{\delta'}^ C(gx))>0$, or that there
exists $h \in \overline {B_{\delta'} ^ C}$ with $hx \in
\overline{B_{\delta'}^ C(gx)}$. In the latter case
$B_{\delta'}^ C(gx)\subset B_{3 \delta'}^ C(x)$. The sets $B_{\delta'}^ C(g)$ for $g \in
B_{2r}^ C$ cover the compact set $\overline{B_{2r}^ C}$. Taking a
finite subcover, we find some $\eta>0$ so that $d(gx,x)>\eta$ or
$gx \in B_{3 \delta'}^ C(x)$ for every $g \in B_{2r}^ C$. It follows
that \eqref{e: intersection equation} holds with $\delta = \min
(\eta, \delta')$. In other words,
$C_r=\bigcup_{\delta>0}D_\delta$, where
$$
D_\delta=\bigl\{ x \in C_r:B_{2r}^ C(x)\cap B_\delta(x)\subset
B_\delta ^ C(x)\bigr\},
$$
and there exists $\delta>0$ with $\mu(D_\delta)>1-\varepsilon$.\smallbreak

Let $K \subset D_\delta$ be compact. We claim that the
$A'$-returns to $K$ are strongly exceptional. So suppose $x \in K$
and $x'=\alpha ^\bs x \in K$ for some $\alpha ^\bs \in A'$. Then
since $x$ and $x '$ are in the same atom of $\mathcal{E}$, the
conditional measures satisfy $\mu_x ^{\mathcal
E}=\mu_{x'}^{\mathcal E}$. By definition of $C_r$ we have
$\mu_x ^{\mathcal E}(B_r ^ C(x))>1/2$ and the same for $x'$.
Therefore $B_r ^ C(x)$ and $B_r ^ C(x')$ cannot be disjoint, and
$x'\in B_{2r}^ C(x)$ follows. By definition of $D_\delta$ it
follows that $x'\in B_\delta ^ C(x)$. Thus the $A'$-returns to $K$
are indeed strongly exceptional.

\demo{$(2) \Longrightarrow (1)$}
Suppose that for every $\ell \geq 1$ there exists a compact set $K_\ell$ with
$\mu(K_\ell)>1-1/\ell$ so that the $A'$-returns to $K$ are strong exceptional.
Then $N=X \setminus \bigcup_\ell K_\ell$ is a nullset. It suffices to show that
(1) holds for every $A'$ ergodic $\mu_x ^{\mathcal E}$ which satisfies
$\mu_x ^{\mathcal E}(N)=0$.

For any such $x$ there exists $\ell>0$ with $\mu_x ^{\mathcal E}(K_\ell)>0$.
Choose some $z \in K_\ell$ with $\mu _ x ^ {\mathcal{E}} (B_{1/m}(z)\cap
K_\ell)>0$ for all $m \geq 1$. We claim that $\mu _ x ^ {\mathcal{E}}$ is
supported on $Cz$, i.e.\ that $\mu _ x ^ \mathcal{E} (Cz) = 1$. Let $\delta$ be
as in the definition of strong exceptional returns. By ergodicity there exists
for $\mu_x ^ {\mathcal{E}}$-a.e.\ $y_0 \in X$ some $\alpha ^\bs \in A'$ with
$y_1=\alpha ^{\bs} y_0 \in B_{\delta}(z)\cap K_\ell$. Moreover, there exists a
sequence $y_n \in A ' y _ 0 \cap K_\ell$ with $y_n \rightarrow z$. Since $y_n
\in B_\delta(y_1)$ for large enough $n$ and since the $A'$-returns to $K_\ell$
are strong exceptional, we conclude that $y_n \in B_\delta ^ C(y_1)$. Since $y_n$
approaches $z$ and $d(z, y _ 1)<\delta$, we have furthermore $z \in \overline
{B_\delta ^ C(y_1)}$. Therefore $y_1 \in Cz$, $y_0=\alpha ^{-\bs}y_1 \in Cz$, and
the claim follows.
\hfq

\begin{lemma}\label{l:Lfoliation}  
{\rm (1)} Under the assumptions of the low entropy case \/{\rm (}\/i.e.\
$s_{ab}(\mu)\break >0$ but $s_{ij}(\mu)=0$ for all $i,j$ with either
$i=a$ or $j=b${\rm ),} there exists a $\mu$-nullset $N \subset X$ such
that for $x \in X \setminus N${\rm ,}
$$
U_{(ab)} x \cap X \setminus N \subset U_{ab} x.
$$

{\rm (2)}
Furthermore{\rm ,} unless $\mu$ is $U_{ab}$-invariant{\rm ,} it can be arranged that
$$
\mu_x ^{ab}\neq \mu_{y}^{ab}
$$
for any $x \in X \setminus N$ and any $y \in U_{(ab)}x \setminus
N$ which is different from $x$.
\end{lemma}

\Proof 
Set $U=U_{(ab)}$ and let $\mu_{x,U}$ be the
conditional measures for the foliation into $U$-orbits. By
\cite[Prop.\ 8.3]{EinsiedlerKatok02} the conditional measure
$\mu_{x,U}$ is a.e.\ -- say for $x\notin N$ --
a product measure of the conditional measures
$\mu_x ^{ij}$ over all $i,j$ for which $U_{ij}\subset U$. Clearly, by the assumptions of the low entropy case,
$\mu_x ^{ab}$ is the only one of these which is nontrivial.
Therefore, $\mu_{x,U}$ -- as a measure on $U$ -- is supported
on the one-dimensional group $U_{ab}$.

By (3) in Section~\ref{sec: conditionalm} the conditional measures satisfy
furthermore that there is a null set -- enlarge $N$ accordingly -- such that
for $x,y\notin N$ and $y=ux \in Ux$ the conditionals $\mu_{x,U}$ and
$\mu_{y,U}$ satisfy that
$\mu_{x,U} \propto \mu _ {y, U} u$.
However, since $\mu_{x,U}$ and $\mu_{y,U}$ are both supported by
$U_{ab}$, it follows that $u \in U_{ab}$. This shows
Lemma~\ref{l:Lfoliation}.(1).

In order to show Lemma~\ref{l:Lfoliation}.(2), we note that we already know
that $y \in U_{ab}x$. So if $\mu_x ^{ab}=\mu_{y}^{ab}$, then $\mu_x ^{ab}$ is
again, by (3) in Section~\ref{sec: conditionalm}, invariant (up to proportionality)
under multiplication by some nontrivial $u \in U_{ab}$. If this were to happen
on a set of positive measure, then by (9) in Section~\ref{sec: conditionalm},
$\mu_x ^{ab}$ are in fact Haar a.e.\ -- a contradiction to our assumption.
\hfq

\Subsec{Sketch of proof of Theorem~{\rm \ref{t: low entropy}}}
We assume that the two equivalent conditions in Proposition \ref{p:
equivalence} fail (the first of which is precisely the condition of
Theorem~\ref{t: low entropy} (2)). From this we will deduce
that $\mu$ is $U_{ab}$-invariant 
which is precisely the statement in Theorem \ref{t: low entropy} (1).

For the following  we may assume without loss of generality that $a=1$
and $b=2$. Write $A'$ and $u(r)=\Id+rE_{12}\in U_{12}$ for $r \in
\RR$ instead of $A_{12}'$ and $u_{12}(r)$. Also, we shall at times
implicitly identify $\mu _ x ^ {1 2}$ (which is a measure on $U _
{12}$) with its push forward under the map $u(r) \mapsto r$, e.g.\ write $\mu _ x ^ {12} ([a,b])$ instead of $\mu _ x ^ {12}
(u([a,b]))$.

By Poincar\'e recurrence we have for a.e.\ $x \in X$ and
every $\delta>0$ that
$$
d(\alpha ^{\bs}x,x)<\delta\mbox{ for some large }\alpha ^{\bs}\in A'.
$$
For a small enough $\delta$ there exists a unique $g \in B_ \delta ^{G}$ such
that $x'=\alpha ^\bs=gx$.

Since $\alpha ^\bs$ preserves the measure and since $A'\subset
L_{12}=C(U_{12})$ the conditional measures satisfy
\begin{equation}\label{e: sketch conditionals}
\mu_x ^{12}=\mu_{x'}^{12}.
\end{equation}
by (5) in Section~\ref{sec: conditionalm}. Since $\mu_x ^{12}$ is nontrivial, we can
find many $r \in \RR$ so that $x(r)=u(r) x$ and $x'(r)=u(r) x'$ are
again typical. By (3) in Section~\ref{sec: conditionalm} the conditionals satisfy
\begin{equation}\label{e: translation}
\mu_{x(r)}^{12} u (r) \propto  \mu_x ^{12}
\end{equation}
and similarly for $x'(r)$ and $x'$. Together with \eqref{e: sketch
conditionals} and the way we have normalized the conditional measures this implies that
$$
\mu_{x(r)}^{12}=\mu_{x'(r)}^{12}.
$$

The key to the low entropy argument, and this is also the key to Ratner's
seminal work on rigidity of unipotent flows, is how the unipotent orbits $x(r)$
and $x'(r)$ diverge for $r$ large (see Figure \ref{Hpfigure}). Ratner's
H-property (which was introduced and used in her earlier works on rigidity of
unipotent flows \cite{Ratner-factors}, \cite{Ratner-joinings} and was generalized by D.
Morris-Witte in \cite{Witte-rigidity}) says that this divergence occurs only
gradually and in prescribed directions. We remark that in addition to our use
of the H-property, the general outline of our argument for the low entropy case
is also quite similar to \cite{Ratner-factors}, \cite{Ratner-joinings}.
\setcounter{figure}{1}
\begin{figure}[t]
\begin{center}
\epsfig{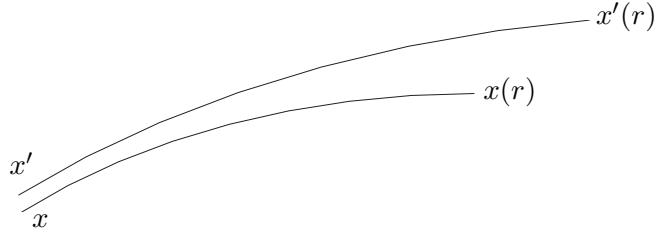}
\setlength{\unitlength}{1mm}
\begin{picture}(0,0)
\put(-80,0){ \put(3,-2){$x$} \put(0,5){$x'$}
\put(63,15){$x(r)$}\put(78,25){$x'(r)$}
}\end{picture}
\end{center}
\caption{\label{Hpfigure} Ratner's H-property: When moving along the unipotent $u(r)$,
the points $x(r)$ and $x'(r)$ noticably differ first only along $U_{(12)}$.}
\end{figure}

We shall use a variant of this H-property in our paper, which at its heart is the following simple matrix calculation (cf.
\cite[Lemma  2.1]{ Ratner-factors} and \cite[Def.~1]{ Ratner-joinings}). Let the entries of $g \in B_\delta ^
G$ be labelled as in
\eqref{e:matrixA}. A simple calculation shows that
$x'(r)=g(r)x(r)$ for
\begin{eqnarray}\label{e:As}
g(r)&=&u(r)gu(-r)\\
&=&\left(
{\begin{array}{ccc}
a_1 + g_{21}r & g_{12}+(a_2-a_1)r -g_{21}r ^ 2 & g_{1*} + g_{2*}r\\
g_{21} & a_2 - g_{21}r & g_{2*}\\
g_{*1} & g_{*2}-g_{*1}r & a_*
\end{array}}
\right).\nn
\end{eqnarray}
Since the return is not exceptional, $g\notin L_{12}=C(U_{12})$
and one of the following holds; $a_2-a_1 \neq 0$, $g_{21}\neq 0$,
$g_{*1}\neq 0$, or $g_{2*}\neq 0$. From this it is immediate that
there exists some $r$ so that $g(r)$ is close to $\Id$ in all
entries except at least one entry corresponding to the subgroup
$U_{(12)}$. More precisely, there is an absolute constant $C$ so
that there exists $r$ with
\begin{gather}\label{e:boundC}
C ^{-1}\leq \max(|(a_2-a_1)r-g_{21}r ^ 2|,\| g_{2*}r \|,\| g_{*1}r
\|)\leq
C,\\
| g_{21}r|\leq C \delta ^{3/8}. \label{e:boundC2}
\end{gather}
With some care we will arrange it so that $x(r),x'(r)$ belong to a fixed
compact set $X_1 \subset X \setminus N$. Here $N$ is as in Lemma
\ref{l:Lfoliation} and $X_1$ satisfies that $\mu_z ^{12}$ depends continuously
on $z \in X_1$, which is possible by Luzin's theorem.

If we can indeed find for every $\delta>0$ two such points
$x(r),x'(r)$ with \eqref{e:boundC} and \eqref{e:boundC2}, we let
$\delta$ go to zero and conclude from compactness that there are
two different points $y,y'\in X_1$ with $y'\in U_{(12)}y$ which
are limits of a sequence of points $x(r),x'(r)\in X_1$. By
continuity of $\mu_z ^{12}$ on $X_1$ we get that
$\mu_y ^{12}=\mu_{y'}^{12}$. However, this contradicts Lemma
\ref{l:Lfoliation} unless $\mu$ is invariant under $U_{12}$.

The main difficulty consists in ensuring that $x(r),x'(r)$ belong
to the compact set $X_1$ and satisfy
\eqref{e:boundC} and \eqref{e:boundC2}. For this we will need several
other compact sets with large measure and various properties.

Our proof follows closely the methods of \cite[\S8]{Lindenstrauss03}.
The
arguments can be simplified if one assumes additional regularity for the
conditional measures $\mu _ z ^{12}$
--- see \cite[\S8.1]{Lindenstrauss03} for more details.

\Subsec{The construction of a nullset and three compact sets}
As mentioned before we will work with two main assumptions: that $\mu$ satisfies the assumptions of the
low entropy case and that the equivalent conditions in Proposition~\ref{p: equivalence} fail. By the former
there exists a nullset
$N$ so that all statements of Lemma~\ref{l:Lfoliation} are
satisfied for $x \in X \setminus N$. By the latter we can assume
that for small enough $\varepsilon$ and for any compact set with
$\mu(K)>1-\varepsilon$ the $A'$-returns to $K$ are not strong
exceptional.

We enlarge $N$ so that $X \setminus N \subset X'$ where $X'$ is as in Section~\ref{sec: conditional}.1. Furthermore, we can assume that $N$ also satisfies
Lemma \ref{l:exceptional}. This shows that for every compact set $K \subset X
\setminus N$ with $\mu(K)>1-\varepsilon$ the $A'$-returns (which exist due to
Poincar\'e recurrence) are not exceptional, i.e.\ for every $\delta>0$ there
exists $z \in K$ and $\bs \in A'$ with $z'=\alpha ^\bs z \in
B_\delta(z)\setminus B_\delta ^ L(z)$.

\demo{Construction of $X_1$} The map $x \mapsto \mu _ x ^ {12}$ is a measurable map from $X$ to a
separable metric space. By Luzin's theorem \cite[p.\ 76]{Fed69} there exists a compact $X_1 \subset X
\setminus N$ with measure
$\mu(X_1)>1-\varepsilon ^ 4$, and the property that $\mu_x ^{12}$ depends
continuously on $x \in X_1$.

\demo{Construction of $X_2$} To construct this set, we use the maximal
inequality (10) in Section~\ref{sec: conditionalm} from \cite[App.~A]{Lindenstrauss03}. Therefore, there exists a
set $X_2 \subset X \setminus N$ of measure $\mu(X_2)>1-C_1 \varepsilon ^ 2$ (with $C_1$ some absolute
constant) so that for any $R>0$ and $x \in X_2$
\begin{equation}\label{e:defX2}
\int_{[-R,R]}1_{X_1}(u(r)x)\dm \mu_x ^{12}(r)\geq
(1-\varepsilon ^ 2)\mu_x ^{12}([-R,R]).
\end{equation}

\demo{Construction of $K=X_3$} Since $\mu_x ^{12}$ is assumed to be
nontrivial a.e., we have $\mu_x ^{12}(\{ 0 \})=0$ and
$\mu_x ^{12}([-1,1])=1$. Therefore, we can find $\rho \in(0,1/2)$ so
that
\begin{equation}\label{e: def X rho}
\mathcal X(\rho)=\bigl\{ x \in X \setminus
N:\mu_x ^{12}([-\rho,\rho])<1/2\bigr\}
\end{equation}
has measure $\mu(\mathcal X(\rho))>1-\varepsilon ^ 2$. Let
$\bt=(1,-1,0,\dots,0)\in \Sigma$ be fixed for the following. By
the (standard) maximal inequality we have that there exists
a compact set $X_3 \subset X \setminus N$ of measure $\mu(X_3)>1-C_2
\varepsilon$ so that for every $x \in X_3$ and $T>0$ we have
\begin{equation}\label{e:defX3}
\begin{split}
\frac{1}{T}\int_0 ^ T1_{X_2}(\alpha ^{-\tau\bt}x)\dm \tau&\geq
(1-\varepsilon),\\
\frac{1}{T}\int_0 ^ T1_{\mathcal
X(\rho)}(\alpha ^{-\tau\bt}x)\dm \tau&\geq
(1-\varepsilon).
\end{split}
\end{equation}

\Subsec{The construction of $z,z'\in X_3${\rm ,} $x,x'\in X_2$}
Let $\delta>0$ be very small (later $\delta$ will approach zero). In
particular, the matrix $g \in B_\delta ^ G$ (with entries as in
\eqref{e:matrixA}) is uniquely defined by $z'=gz$ whenever $z,z'\in X_3$ and
$d(z,z')<\delta$. Since the $A'$-returns to $X_3$ are not exceptional, we can
find $z \in X_3$ and $\alpha ^\bs \in A'$ with $z'=\alpha ^{\bs}z \in B_
\delta(z)\cap X_3$ so that
\begin{equation}\label{e:defkappa_z}
\kappa(z,z')=\max\bigl(|a_2-a_1|,|g_{21}|^{1/2},\| g_{*1}\|,\|
g_{2*}\|\bigr)\in (0,c \delta ^{1/2}),
\end{equation}
where $c$ is an absolute constant allowing us to change from the
metric $d(\cdot,\cdot)$ to the norms we used above.

For the moment let $x=z$, $x'=z'$, and $r=\kappa(z,z')^{-1}$. Obviously
$\max\bigl(|(a_2-a_1)r|,|g_{21}|^{1/2}r,\| g_{2*}r \|,\| g_{*1}r \|\bigr)=1$.
If the maximum is achieved in one of the last two expressions, then
\eqref{e:boundC} and \eqref{e:boundC2} are immediate with $C=1$. However, if the
maximum is achieved in either of the first two expressions, it is possible that
$(a_2-a_1)r-g_{21}r ^ 2$ is very small. In this case we could set $r=2
\kappa ^{-1}(z,z')$, then $(a_2-a_1)r$ is about $2$ and $g_{21}r ^ 2$ is about
$4$. Now \eqref{e:boundC}--\eqref{e:boundC2} hold with $C=10$. The problem with
this naive approach is that we do not have any control on the position of
$x(r),x'(r)$. For all we know these points could belong to the null set $N$
constructed in the last section.

\begin{figure}[ht]
\begin{center}
\epsfig{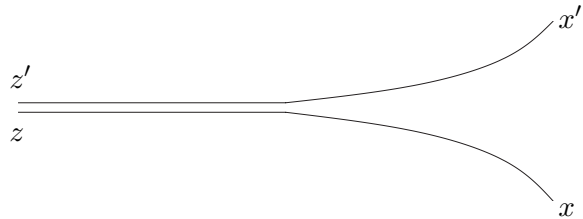}
\setlength{\unitlength}{1mm}
\begin{picture}(0,0)
\put(-90,0){ \put(15,8){$z$}\put(15,15){$z'$} \put(88,-2){$x$}\put(88,23){$x'$}
}\end{picture}
\end{center}
\vglue-12pt
\caption{\label{TrumpFigure} The distance function $\kappa(x,x')$
might be constant for small $\tau$ and increase exponentially
later.} \vglue-12pt
\end{figure}

To overcome this problem we want to use the conditional measure
$\mu_x ^{12}$ to find a working choice of $r$ in some interval
$I$ containing $\kappa(z,z')^{-1}$. Again, this is not immediately
possible since {\it a~priori} this interval could have very small
$\mu_x ^{12}$-measure, or even be a nullset. To fix this, we use
$\bt=(1,-1,0,\dots,0)$ and the flow along the
$\alpha ^{\bt}$-direction in Lemma \ref{l:findtau}. However, note that
$x=\alpha ^{-\tau\bt}z$ and $x'=\alpha ^{-\tau\bt}z'$ differ by
$\alpha ^{-\tau\bt} g \alpha ^{\tau\bt}$. This results possibly in a
difference of $\kappa(x,x')$ and $\kappa(z,z')$ as in Figure \ref{TrumpFigure},
and so we might
have to adjust our interval along the way. The way $\kappa(x,x')$
changes for various values of $\tau$ depends on which terms give
the maximum.

\begin{lemma}\label{ll:findtau}
For $z,z'\in X_3$ as above let $T=\frac{1}{4}|\ln \kappa(z,z')|${\rm ,}
$\eta \in \{ 0,1 \}${\rm ,} and $\theta \in [4T, 6T]$. There exist 
subsets $P,P'\subset[0,T]$ of density at least $1- 9 \varepsilon$
such that for any $\tau \in P$ {\rm (}$\tau \in P'${\rm ),} 
\begin{enumerate}
\item $x=\alpha ^{-\tau\bt}z \in X_2$ ($x'=\alpha ^{-\tau\bt}z'\in
X_2$) and
\item the conditional measure $\mu_x ^{12}$ satisfies the estimate \end{enumerate}
\vglue-20pt
\begin{equation}\label{e:introduceS}
\mu_x ^{12}\bigl([-\rho S,\rho S]\bigr)< \frac{1}{2}\mu_x ^{12}\bigl([-S,S]\bigr)
\end{equation}
\begin{enumerate} \item[]
  where $S=S(\tau)=e ^{\theta-\eta \tau}$ \/{\rm (}\/and similarly for
$\mu_{x'}^{12}${\rm )}.
\end{enumerate}
\end{lemma}

\Proof 
By the first line in \eqref{e:defX3} there exists a set $Q_1
\subset [0,T]$ of density at least $1-\varepsilon$ (with respect to
the Lebesgue measure) such that $x=\alpha ^{-\tau\bt}z$ belongs to
$X_2$ for every $\tau \in Q_1$.

By the second line in \eqref{e:defX3} there exists a set $Q_2
\subset[0, 4T]$ of density at least $1-\varepsilon$ such that
$\alpha ^{-v\bt} z \in \mathcal X(\rho)$ for $v \in Q_2$. Let
$$
Q_3=\Bigl\{ \tau \in[0,T]:\frac{1}{2}(\theta+(2-\eta)\tau) \in
Q_2\Bigr\}.
$$
A direct calculations shows that $Q_3$ has density at least $1- 8
\varepsilon$ in $[0,T]$, and for $\tau \in Q_3$ and
$v=\frac{1}{2}(\theta+(2-\eta)\tau)$ we have $y=\alpha ^{-v\bt}z
\in \mathcal X_\rho$.

\begin{figure}[t]
\vglue12pt
\begin{center}
\epsfig{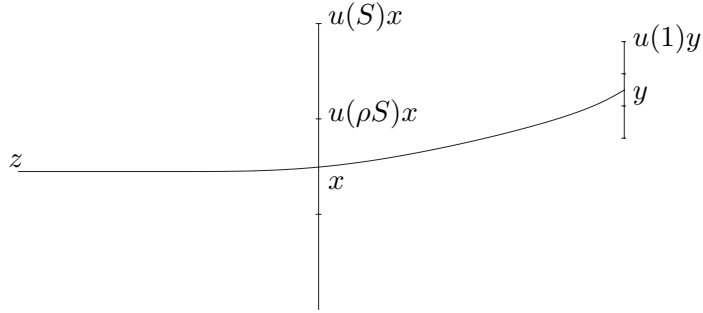}
\setlength{\unitlength}{1mm}
\begin{picture}(0,0)
\put(-60,0){ \put(17.5,16){$x$} \put(-25,19){$z$}\put(58,28){$y$}
\put(17.5,38){$u(S)x$}\put(17.5,25){$u(\rho S)x$}
\put(58,35){$u(1)y$}
}\end{picture}
\end{center}
\caption{\label{stretchfig} From the way the leaf $U_{12}x$ is contracted along
$\alpha ^{-\bt}$ we can ensure \eqref{e:introduceS} if $y=\alpha ^{-w\bt}x \in
\mathcal X_\rho$} \vglue-12pt
\end{figure}
We claim the set $P=Q_1 \cap Q_3 \subset[0,T]$ satisfies all
assertions of the lemma; see Figure \ref{stretchfig}. First $P$ has at least density $1-9
\varepsilon$. Now suppose $\tau \in P$; then $x=\alpha ^{-\tau\bt}z
\in X_2$ by definition of $Q_1$. Let $w=\frac{1}{2}(\theta-\eta
\tau)$; then
$$
y=\alpha ^{-w\bt}x=
\alpha ^{-v\bt}z \in \mathcal X_\rho
$$
by the last paragraph. By \eqref{e: def X rho}
$$
\mu_y ^{12}([-\rho,\rho])<\frac{1}{2}\mu_y ^{12}([-1,1])=\frac{1}{2}.
$$
By property (4) in  Section \ref{sec: conditional} of the conditional
measures we get that
$$
\frac {\mu _ y ^ {12} ([- \rho, \rho]) }{ \mu _ y ^ {1 2} ([-
1,1])} = \frac {(\alpha ^ {-w \mathbf t} \mu _ x ^ {12} \alpha ^
{w \mathbf t})([- \rho, \rho]) }{ (\alpha ^ {-w \mathbf t} \mu _ x
^ {12} \alpha ^ {w \mathbf t})([- 1, 1])} = \frac {\mu _ x ^ {12}
([- \rho e ^ {2w}, \rho e ^ {2w}]) }{ \mu _ x ^ {1 2} ([- e ^
{2w},e ^ {2w}])}
$$
This implies \eqref{e:introduceS} for
$S=e ^{2w}=e ^{\theta-\eta \tau}$. The construction of $P'$ for $z'$
is similar.
\Endproof\vskip4pt  

The next lemma uses Lemma~\ref{ll:findtau} to construct $x$ and $x'$ with the
property that certain intervals containing $\kappa(x,x')^{-1}$
have $\mu_x ^{12}$-measure which is not too small. This will allow us in Section~\ref{ss: conclusion} to
find $r$ so that both $x(r)$ and $x'(r)$ have all the desired properties.

\begin{lemma}\label{l:findtau}
Let $z, z ' \in X _ 3$ and $T = \frac {1 }{ 4} \absolute {\ln \kappa (z, z ')}$
be as above. If $\varepsilon<\frac{1}{100}${\rm ,} then there exists $\tau \in [0,T]$
such that
\begin{enumerate}
\item both $x=\alpha ^{-\tau\bt}z$ and
$x'=\alpha ^{-\tau\bt}z'$ are in $X_2${\rm ,}
\item
$\kappa(x,x')<c \delta ^{3/8}${\rm ,} and
\item for \( R=\kappa(x,x')^{-1}
\) \/{\rm (}\/as well as $R'=\rho ^{-5}R$\/{\rm )}\/ 
\end{enumerate} 
\vglue-20pt
\begin{equation}\label{e:introduceR}
\begin{aligned}
\mu_x ^{12}\bigl([-\rho R,\rho R]\bigr)&<\frac{1}{2}\mu_x ^{12}\bigl([-R,R]\bigr)\mbox{ and }\\
\mu_{x'}^{12}\bigl([-\rho R,\rho
R]\bigr)&<\frac{1}{2}\mu_{x'}^{12}\bigl([-R,R]\bigr).
\end{aligned}
\end{equation}
\end{lemma}

\Proof 
Let
\begin{eqnarray*}
\kappa_a(z,z')&=&|a_2-a_1|,\\
\kappa_u(z,z')&=&\max\bigl(|g_{21}|^{1/2},\|g_{*1}\|,\|g_{2*}\|\bigr)\in
(0,c\delta^{1/2}).
\end{eqnarray*}
The corresponding quantities for $x,x'$ are defined similarly. The
number $T$ is chosen so that the two points $x=\alpha ^{-\tau\bt}z$
and $x'=\alpha ^{-\tau\bt}z'$ are still close together for $\tau
\in[0,T]$. In fact,
\begin{equation}\label{e: def tilde g}
\tilde{g}= \alpha ^{-\tau\bt} \left(
{\begin{array}{ccc} a_1 & g_{12} & g_{1*}\\
g_{21} & a_2 & g_{2*}  \\
g_{*1} & g_{*2} & a_*
\end{array}}
\right)\alpha ^{\tau\bt}= \left(
{\begin{array}{ccc} a_1 & e ^{-2 \tau}g_{12} & e ^{-\tau}g_{1*}\\
e ^{2 \tau}g_{21} & a_2 & e ^{\tau}g_{2*}  \\
e ^{\tau}g_{*1} & e ^{-\tau}g_{*2} & a_*
\end{array}}
\right)
\end{equation}
satisfies $x'=\tilde{g} x$, and so
\begin{eqnarray}\label{e:findkappa_x} \quad 
\kappa_a(x,x')&=&\kappa_a(z,z'),\
\kappa_u(x,x')=e ^ \tau \kappa_u(z,z')\leq \kappa(z,z')^{\frac{3}{4}}<c \delta ^{\frac{3}{8}}\mbox{ and }
\\
\kappa(x,x')&=&\max(\kappa_a(x,x'),\kappa_u(x,x'))<c \delta ^{\frac{3}{8}};\nn
\end{eqnarray}
see also Figure \ref{TrumpFigure}. Hence the second statement of the lemma holds.

For the other two statements of the lemma we will use Lemma \ref{ll:findtau} to
define four subsets $P _ a, P _ u, P ' _ a, P ' _ u \subset [0,T]$, each of
density at least $1 - 9 \varepsilon$, so that for every $\tau$ in the intersection
of these four sets both (1) and (3) hold.

\demo{Definition of $P_a$} If $\kappa (x,x ') > \kappa _ a (x,x
')$ for all $\tau \in [0,T]$ (recall that $x,x'$ depend implicitly
on $\tau$) we set $P _ a = [0, T]$.

Otherwise, it follows from \eqref{e:findkappa_x} that $\kappa (z, z ') = \kappa _ a (z, z ')$.
We apply Lemma \ref{ll:findtau} for $\eta=0$ and $\theta=-\log
\kappa_a(z,z')=4T$, and see that \eqref{e:introduceS} holds for $\tau \in P_a$, where $P_a \subset
[0,T]$ has density at least $1- 9 \varepsilon$, and
$S_a=\kappa_a(z,z')^{-1} = \kappa _ a (x,x ') ^{-1}$.

\demo{Definition of $P _ u$} If $\kappa (x,x ') > \kappa _ u (x,x ')$ for all $\tau \in [0,T]$ we set $P _ u = [0,
T]$.

Otherwise, it follows from \eqref{e:findkappa_x} that
$\kappa(x,x')=e ^ \tau \kappa_u(z,z')\geq \kappa(z,z')$ for some
$\tau \in[0,T]$; hence $\kappa_u(z,z')\in
[\kappa(z,z')^{{5}/{4}},\kappa(z,z')]$. This time, we apply Lemma
\ref{ll:findtau} with $\eta=1$ and $\theta=-\log \kappa_u(z,z')\in
[4T,5T]$. We conclude that in this case \eqref{e:introduceS} holds
for $\tau \in P_u$, where $P_u \subset[0,T]$ is a set of density
$1- 9 \varepsilon$, and $ S_u=\kappa_u(x,x')^{-1}=e ^{\theta-\tau}$.

Clearly, since $\kappa (x,x ')$ is either $\kappa _ u (x,x ')$ or $\kappa _ a
(x,x ')$ at least one of the sets $P_a$ or $P_u$ is constructed using Lemma
\ref{ll:findtau}; so in particular if $\tau \in P _ a \cap P _ u$ then $x \in
X_2$. Furthermore, if $\tau \in P _ a \cap P _ u$ we have that
\eqref{e:introduceS} holds for $S=R=\kappa(x,x')^{-1} = \min (\kappa_a (x,x '),
\kappa _ u(x, x'))$.

The sets $P_a'$ and $P_u'$ are defined similarly using $z'$.

The set $P_a \cap P_a'\cap P_u \cap P_u'\subset [0,T]$ has
density at least $(1- 36 \varepsilon)$, so that  in particular if $\varepsilon$ is small it is nonempty.
For any $\tau$ in this intersection, $x,x'\in X_2$ and \eqref{e:introduceR} holds for $R=\kappa(x,x')^{-1}$.

The additional statement in the parenthesis follows similarly,
the only difference being the use of a slightly different value for $\theta$ in both
cases, and then taking the intersection of $P_a \cap P_a'\cap P_u \cap P_u'$ with four more subsets of
$[0,T]$ with similar estimates on their densities.
\hfq

\Subsec{Construction of $x(r)${\rm ,} $x'(r)$ and the conclusion of the proof} \label{ss: conclusion}
Recall that we found $z,z'\in X_3$ using Poincar\'e recurrence and the
assumption that the $A'$-returns to $X_3$ are not exceptional. In the last
section we constructed $x=\alpha ^{-\tau\bt}z,x'=\alpha ^{-\tau\bt}z'=\alpha ^\bs
\in X_2$ using the properties of $X_3$ to ensure \eqref{e:introduceR}. Since
$\alpha ^{\bs}$ acts isometrically on the $U_{12}$-leaves, it follows from
property (4) of the conditional measures in Section~\ref{sec: conditional} that
$\mu_z ^{12}=\mu_{z'}^{12}$ and $\mu_x ^{12}=\mu_{x'}^{12}$.

Let
$$
\begin{aligned}
P&=\{ r \in [-R,R]: u(r)x \in X_1 \}\mbox{ and}\\
P'&=\{ r \in [-R,R]: u(r)x'\in X_1 \}.
\end{aligned}
$$
By \eqref{e:defX2} we know that $P$ and $P'$ both have density at
least $(1- \varepsilon ^ 2)$ with respect to the measure
$\mu_x ^{12}=\mu_{x'}^{12}$. By \eqref{e:introduceR} we know that
$[-\rho R,\rho R]$ contains less than one half of the
$\mu_x ^{12}$-mass of $[-R,R]$. Therefore, if $\varepsilon$ is small
enough there exists $r \in P \cap P'\setminus[-\rho R,\rho R]$. We
define $x(r)=u(r)x$ and $x'(r)=u(r)x'$, and conclude that
$x(r),x'(r)\in X_1$ satisfy $\mu_{x(r)}^{12}=\mu_{x'(r)}^{12}$ by
property (3) in Section~\ref{sec: conditional}.

Let $\tilde g$ be defined as in \eqref{e: def tilde g} and write
$\tilde g_{12},$ \dots for the matrix entries. With $\tilde
g(r)=u(r)\tilde gu(-r)$ we have $x'(r)=\tilde g(r) x(r)$ and
$$
\tilde g(r)=\left( {\begin{array}{ccc}
a_1 + \tilde g_{21}r & \tilde g_{12}+(a_2-a_1)r -\tilde g_{21}r ^ 2 & \tilde g_{1*} + \tilde g_{2*}r\\
\tilde g_{21} & a_2 - \tilde g_{21}r & \tilde g_{2*}\\
\tilde g_{*1} & \tilde g_{*2}-\tilde g_{*1}r & a_*
\end{array}}
\right).
$$
We claim it is possible to achieve
\begin{gather}\label{e: bound C again}
C ^{-1}\leq \max(|(a_2-a_1)r-\tilde g_{21}r ^ 2|,\| \tilde g_{2*}r
\|,\| \tilde g_{*1}r \|)\leq
C,\\
| \tilde g_{21}r| \leq C \delta ^{3/8}. \label{e: bound C2 again}
\end{gather}
for some constant $C$; see Figure \ref{Hpfigure}.

We proceed to the proof of \eqref{e: bound C again} and \eqref{e: bound C2 again}.
For \eqref{e: bound C2 again} we first recall that $|\tilde
g_{21}|\leq \kappa(x,x')^ 2$, and then use \eqref{e:defkappa_z} and
\eqref{e:findkappa_x} to get
$$
|\tilde g_{21}r|\leq \kappa(x,x')^ 2R=\kappa(x,x')\leq
\kappa(z,z')^{3/4}\leq (c \delta ^{1/2})^{3/4} \leq
c ^{3/4}\delta ^{3/8}.
$$

We now turn to prove \eqref{e: bound C again}. It is immediate
from the definition of $R$ that
\begin{equation}\label{e:boundwrho}
\rho \leq \max\bigl(|(a_2-a_1)r|,|\tilde g_{21}|^{1/2}|r|,\|
\tilde g_{2*}r \|,\| \tilde g_{*1}r \|\bigr)\leq 1.
\end{equation}
There are two differences of this estimate to the one in \eqref{e:
bound C again}; first we need to take the square of the second
term -- this replaces the lower bound $\rho$ by its square,
secondly we looked above at $(a_2-a_1)r$ and $\tilde g_{21}r ^ 2$
separately -- taking the difference as in \eqref{e:boundC} might
produce a too small a number (almost cancellation). So \eqref{e:
bound C again} follows with $C=\rho ^{-2}$, unless
\begin{gather} \label{e: first bad case}
\max\bigl(|(a_2-a_1)r|,|\tilde g_{21}|^{1/2}|r|\bigr)\geq \rho,\\
|(a_2-a_1)r-\tilde g_{21}r ^ 2|<\rho ^ 2<\frac{\rho}{2}. \label{e:
second bad case}
\end{gather}

This is a minor problem, and we can overcome it using the last
statement in Lemma \ref{l:findtau}. Assume that for some $r \in
[-R,R]\setminus[-\rho R,\rho R]$ this problem occurs. We deduce a
lower estimate on $|\tilde g_{21}|$. If the maximum in \eqref{e:
first bad case} is achieved at $|\tilde g_{21}|^{1/2}|r|\geq
\rho$, then $|\tilde g_{21}|\geq \rho ^ 2 R ^{-2}$. If the maximum is
achieved at $|(a_2-a_1)r|\geq \rho$, \eqref{e: second bad case}
shows that $|\tilde g_{21}r ^ 2|\geq \rho/2 \geq \rho ^ 2$ (since
$\rho<1/2$) and so in both cases
\begin{equation}\label{e:badcase2}
|\tilde g_{21}|\geq \rho ^ 2 R ^{-2}.
\end{equation}
Now we go through the construction of $r$ again,
only this time using the last statement in Lemma \ref{l:findtau},
and find
$r \in[\rho ^{-5}R,\rho ^{-5}R]\setminus[\rho ^{-4}R,\rho ^{-4}R]$. The
equivalent to \eqref{e:boundwrho} is now the estimate
$$
\rho ^{-4}\leq \max\bigl(|(a_2-a_1)r|,|\tilde g_{21}|^{1/2}|r|,\|
\tilde g_{2*}r \|,\| \tilde g_{*1}r \|\bigr)\leq \rho ^{-5}.
$$
This shows that $|(a_2-a_1)r|\leq \rho ^{-5}$, and
\eqref{e:badcase2} shows that
$$
|\tilde g_{21}r ^ 2|\geq \rho ^ 2 R ^{-2} (\rho ^{-4}R)^ 2=\rho ^{-6}.
$$
Together, we find a lower bound for
$$
|(a_2-a_1)r-\tilde g_{21}r ^ 2|\geq \rho ^{-6}-\rho ^{-5}>0,
$$
i.e.\ the problem of almost cancellation cannot happen again.

Starting with the nonexceptional return $z'\in X_3$ of $z \in X_3$
we have found two points $x_r,x_r'\in X_1$  which satisfy
\eqref{e: bound C again}, \eqref{e: bound C2 again}. Since we assume to have
nonexceptional returns to $X_3$ for every $\delta=\frac{1}{n}>0$,
we get two sequences $y_n$ and $y_n'$ of points in $X_1$ with the
same conditional measures
$$
\mu_{y_n}^{12}=\mu_{y_n'}^{12}.
$$
Compactness shows that we can find convergent subsequences with
limits $y,y'\in X_1$. It follows from \eqref{e: bound C
again}, \eqref{e: bound C2 again} that $y'\in Uy$, and from
\eqref{e: bound C again} that $y'\neq y$. By continuity of
$\mu_x ^{12}$ for $x \in X_1$ the conditional measures
$\mu_y ^{12}=\mu_{y'}^{12}$ agree. However, this contradicts Lemma
\ref{l:Lfoliation}, unless $\mu$ is invariant under $U_{12}$.

\section{Proof that exceptional returns are not possible for
$\Gamma = \SL (k, \ZZ)$} 
\label{s: exceptional returns}

If case (3) in Theorem~\ref{t:slk} holds, then this gives some
restriction on $\Gamma$. In other words, for some lattices in $\SL
(k, \RR)$, exceptional returns cannot occur. As will be shown
below, such is the case for $\Gamma = \SL (k, \ZZ)$.

We recall that $H_{ab} \subset \SL (k, \RR)$ is an $A$-normalized subgroup
isomorphic to $ \SL (2, \RR)$, and $A ' = A \cap \centralizer {H _ {ab}}$. If
case (3) of Theorem~\ref{t:slk} holds then any\break $A '$-ergodic component of $\mu$
is supported on a single $ \centralizer {H _ {ab}}$-orbit. In particular, we
have an abundance of $A '$-invariant probability measures supported on single
$\centralizer {H _ {ab}}$ orbits. The mere existence of such measures is a
restriction on $\Gamma$.

\begin{theorem}
\label{theorem about exceptional returns}
Suppose that $ \nu$ is an $A'$ invariant probability measure on $X = \SL (k, \RR) / \Gamma${\rm ,}
 and that $\supp \nu \subset \centralizer {H _ { ab }}x$ for some $ x \in X$. Then there is a $ \gamma \in
\Gamma$ which is
\begin{enumerate}
\item diagonalizable over $\RR,$
\item $\pm 1$ is not an eigenvalue of $\gamma.$
\item All eigenvalues of $\gamma$ are simple except precisely one which has multiplicity two.
\end{enumerate}
\end{theorem}

Before we prove this theorem, we note the following:

\begin{proposition}\label{proposition about special linear group}
There is no  $ \gamma \in \SL (k, \ZZ)$ satisfying the three conditions of
 Theorem~{\rm \ref{theorem about exceptional returns}}.
\end{proposition}

In particular, case (3) of Theorem~\ref{t:slk} cannot occur for $\SL (k, \RR) / \SL (k, \ZZ)$.

\demo{Proof of Proposition~{\rm \ref{proposition about special linear group}}}
Suppose $\gamma \in \SL (k, \ZZ)$ is diagonalizable over $\RR$. Then its eigenvalues (with the correct multiplicities) are roots of the characteristic polynomial of $\gamma$, a polynomial with integer coefficients and both leading term and constant term equal to one. If there is some eigenvalue which is not equal to $\pm 1$ and which occurs with multiplicity greater than one then necessarily this eigenvalue is not rational, and its Galois conjugates would also have multiplicity greater than one, contradicting (3).
\Endproof\vskip4pt

To prove Theorem~\ref{theorem about exceptional returns}, we need the following standard estimate:

\begin{lemma}
\label{a lemma about eigenvalues} There is a neighborhood $U_0$ of
the identity in $\SL (m, \RR)$ so that for any $\lambda _ 1,
\lambda _ 2, \dots, \lambda _ m$ with $ \absolute {\lambda _ i -
\lambda _ j} > 1$ and $h \in U_0$ one has that\break $ h \diag (e ^
{\lambda _ 1}, e ^ {\lambda _ 2}, \dots, e ^ {\lambda _ m})$ is
diagonalizable over $\RR$ with positive eigenvalues and the
eigenvalues $ e ^ {\lambda ' _ 1}, e ^ {\lambda ' _ 2}, \dots, e ^
{\lambda ' _ m}$ satisfy $\absolute {\lambda ' _ i - \lambda _ i}
< \tfrac{1}{2}$.
\end{lemma}

\Proof 
Without loss of generality, suppose $\lambda _ 1 > \lambda _ 2 >
\dots > \lambda _ m$. Let $f=\diag (e ^ \lambda _ 1, e ^ \lambda _
2, \dots, e ^ \lambda _ m)$  and $\eta _ 1, \eta _ 2, \dots$ be
the eigenvalues of $f' = h f$ ordered according to descending
absolute value. Set for $ 1 \leq i \leq m$, $\lambda ' _ i = \log
\absolute {\eta _ i}$.

Clearly, $\lambda' _ 1 = \lim _ {n \to \infty} \frac {\log \norm
{{f'} ^ n} }{ n}$. Since $f$ is self adjoint, $\norm f = e ^
{\lambda _ 1}$ so
$$
\lambda ' _ 1 \leq \log \norm h + \lambda _ 1 .$$ Let
$ \delta > 0 $ be small (it will be chosen later and will be
independent of $f $). Consider the cones
\begin{align*}
K= & \left\{ (x _ 1, \dots, x _ m): \absolute {x _ l} \leq \delta \absolute {x _ 1}
\text{for every $l \ne 1$} \right\} \\
K ' =& \left\{(x _ 1, \dots, x _ m): \absolute {x _ l} \leq \delta
e ^ {- 1} \absolute {x _ 1} \text{for every $l \ne 1$} \right\}
.\end{align*} Then $f K \subset K ' $, and for every $x \in K $
$$
\norm {fx} \geq (1 - c \delta) e ^ {\lambda _ 1} \norm x
,$$ for some $c$ depending only on $m$. Suppose now that $h$ is close enough to the identity so that $ h K ' \subset K $.
Then $ f ' K = hf K \subset h K ' \subset K $. Again assuming
that $h$ is in some fixed neighborhood of the identity, for any $x
\in K$, we have
$$
\norm {f ' x} \geq (1 - 2 c \delta) e ^ {\lambda _ 1}\norm x
$$
so that
$$
\norm {{f '} ^ n} \geq \norm {{f '} ^ ne _ 1} \geq \bigl((1 - 2 c
\delta) e ^ {\lambda _ 1})\bigr) ^ n .$$ In other
words, if $h$ is in some fixed neighborhood of the identity
(independently of $f$) then
$$
\absolute {\lambda ' _ 1 - \lambda _ 1} < C_1 \delta
.$$

Similarly, $e ^{\lambda _ 1 + \lambda _ 2}$ is the dominating
eigenvalue of $ f \wedge f $, i.e.\ the natural action of $f$ on
the space $ \RR ^ n \wedge \RR ^ n $. Applying the same logic as
before, $\lambda _ 1' + \lambda _ 2' = \lim _ {n \to \infty} \frac
{\log \norm {(f' \wedge f') ^ n} }{ n}$, and as long as $h$ is in
some fixed neighborhood of the identity, independently of $f$
$$
\absolute {\lambda ' _ 1 + \lambda ' _ 2 - \lambda _ 1 - \lambda _ 2} < C _ 2 \delta
$$
and more generally
\begin{equation} \label{lambda prime close to lambda}
\absolute {\sum_ {i = 1} ^ k (\lambda ' _ i - \lambda _ i)} < C _
i \delta .\end{equation} Clearly, \eqref{lambda prime close to
lambda} implies that there is some $C$ depending only on $m$, and
a neighborhood of the identity in $\SL (m, \RR)$ depending only on
$ \delta$ so that if $h$ is in that neighborhood
$$
\absolute {\lambda ' _ i - \lambda _ i} < C \delta \qquad
\text{for all $1 \leq i \leq m$} .$$ In particular,
if $\lambda _ i > \lambda _ {i - 1} + 1$ for every $i$ then if $C
\delta < \tfrac{1}{2}$,  all $ \lambda ' _ i$ are distinct. Since
this holds for all $h$ in a connected neighborhood of the
identity, all the eigenvalues of $ f '$ are real and also
positive, so that $ \eta _ i = e ^ {\lambda ' _ i}$.
\hfq

\demo{Proof of Theorem~{\rm \ref{theorem about exceptional returns}}}
Without loss of generality, we may take $a = 1,\break b = 2$. Let $ a (t) = \diag (e ^
{\lambda _ 1 t}, \dots, e ^ {\lambda _ kt})\in A'$ with $\lambda _ 1 = \lambda
_ 2$ and for every other pair $i, j$ we have $\lambda _ i \neq \lambda _ j$.

Take $U _ 0$ to be a symmetric neighborhood of the identity in
$\SL (k - 2, \RR)$ as in Lemma~\ref{a lemma about eigenvalues},
and
$$
U _ 1 = \left\{
\begin{pmatrix}
e ^ r& 0& 0 \\
0& e ^ r& 0 \\
0& 0& e ^{-\frac{2r}{k-2}}h
\end{pmatrix}
: \text {$r \in (- 1 /8, 1 / 8)$ and $ h \in U_0$} \right\} ,$$
and let $t _ 0 = 2 \max_ {\lambda _ i \neq \lambda _ j}  \absolute {\lambda _ i
- \lambda _ j} ^{-1}$. Note that $U _ 1$ is also symmetric, i.e.\break $U _ 1 ^{-1}
= U _ 1$.

By Poincar\'e recurrence, for $\nu$-almost every $x = g \Gamma \in
\SL (k, \RR) / \Gamma$ there is a $t > t _ 0$ so that $ a (t) x
\in U _ 1 x$; so in particular $U _ 1 a(t) \cap g \Gamma g ^{-1}
\neq \emptyset$. Let
$$
\tilde \gamma = \begin{pmatrix} e ^ s& 0& 0\\
0& e ^ s& 0\\
0& 0& f'
\end{pmatrix} \in U _ 1 a(t) \cap g \Gamma g ^{-1}
$$
be any element from this intersection.
By assumption, for every pair $i, j$ except $1,2$ we have
\begin{equation} \label{inequality about absolute values}
\absolute {\lambda _ i - \lambda _ j} t > 2,
\end{equation}
and we can apply Lemma~\ref{a lemma about eigenvalues} to deduce that the
eigenvalues of
$$
f' = h e ^{-\frac{2s}{k-2}} \diag (e ^ {\lambda _ 3t}, e ^
{\lambda _ 4 t}, \dots, e ^ {\lambda _ kt})
$$
for some $h \in U _ 0$ are of the form $e ^ {\lambda _ 3'}, \dots, e ^ {\lambda
_ k'}$ with $\absolute {\lambda _ i' - \lambda _ it} < 3 / 4$ for
$i=3,\dots,k$. Finally $\absolute {s - \lambda _ 1t} = \absolute {s - \lambda
_ 2 t} \leq 1 / 8$.

In view of \eqref{inequality about absolute values} it is clear
that $\tilde \gamma$ and hence $g ^{-1} \tilde \gamma g \in
\Gamma$ satisfy all the conditions of Theorem~\ref{theorem about
exceptional returns}.
\hfq

\section {Conclusion of the proof of Theorem~\ref{theorem about lattice space}}
In this section, we conclude the derivation of
Theorem~\ref{theorem about lattice space}, and its corollary,
Corollary~\ref{corollary about lattice space}, from
Theorem~\ref{t:slk}.  Throughout this section, $X$ will denote the
quotient space $\SL (k, \RR) / \SL (k, \ZZ)$, and $\mu$ be an $
A$-ergodic and invariant probability measure on $X$. For every
pair $a, b$ of distinct indices in $\left\{ 1, \dots, k \right\}
$, one of the three possibilities of Theorem~\ref{t:slk} holds.
However, in view of the results of the previous section, in
particular Theorem~\ref{theorem about exceptional returns} and
Proposition~\ref{proposition about special linear group},
Theorem~\ref{t:slk}.(3), i.e.\ the case of exceptional returns,
cannot occur for the lattice $\SL (k, \ZZ)$. Therefore, for every
pair $ a, b $ of distinct indices one of the following two
mutually exclusive possibilities holds:
\begin{enumerate}
\item
The conditional measures $\mu _ x ^ {ab}$ and $\mu _ x ^ {ba}$ are trivial a.e.
\item
The conditional measures $\mu _ x ^ {ab}$ and $\mu _ x ^ {ba}$ are Haar and $\mu$ is invariant under left multiplication with elements of $H _ {ab} = \langle U _ {ab}, U _ {ba} \rangle$.
\end{enumerate}

Define a relation $a \sim b$ if $ \mu$ is $U _ {ab}$-invariant. By (2) above it
follows that $a \sim b$ if and only if $b \sim a$. Furthermore, since the group
generated by $U _ {ab}$ and $ U _ {bc}$ contains $ U _ {ac}$, it is clear that
$\sim$ is in fact an equivalence relation on $\left\{ 1, \dots, k \right\}$.
Let $H$ be the group generated by all $U _ {ab}$ with $a \sim b$. Let $r$
denote the number of equivalence classes for $\sim$ which contain more than
one element, and $k _ 1, k _ 2, \dots, k _ r$ be their sizes; so in particular
$\sum_ {i = 1 } ^ rk _ i \leq k$. By permuting the indices if necessary we can
assume these equivalence classes are consecutive indices and $H = \prod_ {i =
1} ^ r \SL (k _ i, \RR)$. By definition, $H$ leaves the measure $\mu$ invariant, is normalized by $A$, and is generated by unipotent one-parameter subgroups of $\SL (n, \RR)$ --- indeed, $H$ is precisely the maximal subgroup of $\SL (n, \RR)$ satisfying these three conditions.

Measures invariant under groups generated by unipotent one-parameter groups are
well understood. In particular, in a seminal series of papers culminating in
\cite{Ratner-measure-rigidity}, M. Ratner showed that if $H$ is such a group
the only $H$-ergodic and invariant probability measures are the algebraic
measures: $L$-invariant measures supported on a closed $L$-orbit for some $L >
H$ (here and throughout, we use the notation $L>H$ to denote that $H$ is a
subgroup and $L \supset H$; specifically $H$ may be equal to $L$) . For the
$A$-invariant measure $\mu$ and $H$ as above we only know that $\mu$ is
$A$-ergodic and $H$-invariant, but not necessarily $H$-ergodic; we shall use
the following version of these measure-rigidity results by Margulis and Tomanov
\cite{Margulis-Tomanov-2}{}\footnote{The main theorem of
\cite{Margulis-Tomanov-2} was substantially more general than what we quote
here. In particular, in their theorem $\Gamma$ can be any closed subgroup of
$G$, and the group $G$ can be a product of real and $p$-adic Lie groups
(satisfying some mild additional conditions).} (similar techniques were used
also in \cite [proof of Thm.~1]{Mozes-epimorphic}; see also
\cite[\S4.4]{Kleinbock-et-al} and \cite{Starkov-strict-ergodicity}). For
any connected real Lie group $G$, we shall say that $g \in G$ in an element of
class $\mathcal{A}$ if $\Ad g$ is semisimple, with all eigenvalues integer
powers of some  $\lambda \in \RR \setminus \left\{ \pm 1 \right\}$, and $g$ is
contained in a maximal reductive subgroup of $G$.

\begin{theorem} [{\cite [Thms.~(a) and (b)] {Margulis-Tomanov-2}}] \label{version of Ratner}
Let $G$ be a connected real Lie\break group{\rm ,} $\Gamma < G$ a discrete subgroup{\rm ,} and
$\tilde H$ generated by unipotent one-parameter groups and elements of class
$\mathcal{A}${\rm ,}
 with $H< \tilde H$ the subgroup generated by uni\-potent one-parameter groups. Let $\mu$
be an $\tilde H$-invariant and ergodic probability measure on $G / \Gamma$. Then there is an $L \geq H$ so
that almost every
$H$-ergodic component of $\mu$ is the $L$-invariant probability measure on a
closed $L$-orbit. Furthermore{\rm ,} if
$$
SN_G(L)= \left\{ g \in N_G(L): \text{conjunction by $g$ preserves Haar measure on $L$}\right\}
$$
then $\tilde H < SN_G(L)$ and $\mu$ is supported on a single $SN_G(L)$-orbit.
In particular{\rm ,} $L$ is normalized by $\tilde H$.
\end{theorem}

\begin{lemma} \label{l: full rank}
Let $H = \prod_ {i = 1} ^ r \SL (k _ i, \RR)$ with $\sum_ {i = 1} ^ rk _ i <
k$. Let
$$
X _ H = \left\{ x \in X: H x \text{ is closed and of finite
volume} \right\} .$$ Then there is a one-parameter
subgroup $a (t)$ of $A$ so that for every $x \in X _ H$ its
trajectory $a (t) x \to \infty$ as $t \to \infty$.
\end{lemma}

\Proof 
Suppose $H x$ is closed and of finite volume, with $x = g \SL (k,
\ZZ)$ and $g=(g _ {ij})$. Let $k ' = \sum_ i k _ i < k$. Let
$\Lambda = g ^{-1} H g \cap \Gamma$. Since $\Lambda$ is Zariski
dense in $g ^{-1} H g$ there is a $\gamma = g ^{-1} h _ 0 g \in
\Lambda$ with $ h _ 0 =
\begin{pmatrix} h ' _ 0& 0 \\ 0& I _ {k - k '} \end{pmatrix}$ so
that \vglue-14pt
\begin{equation} \label{rationality}
V_g := \left\{  y \in \RR ^ k: y ^ T g ^{-1} h g = y ^ T \quad
\text{for all $h \in H$} \right\} = \left\{ y \in \RR ^ k: y ^ T
\gamma = y ^ T \right\} .\end{equation} Notice that since $g (g
^{-1} h _ 0 g) = h _ 0 g$  (the transpose of) the last $k - k '$
rows of $ g$ are in $V _ g$.

Clearly $\dim V _ g = k - k '$, and using the right hand side of
\eqref{rationality} it is clear that $V _ g$ is a rational
subspace of $\RR ^ n$ (i.e.\ has a basis consisting of rational
vectors). Since $V _ g$ is rational, there is an integer vector $m
\in \ZZ ^ n \cap (V _ g) ^ \perp$. In particular, the last $k - k
'$ entries in the vector $gm$ (which is a vector in the lattice in
$\RR ^ k$ corresponding to $g \SL (k, \ZZ)$) are zero.  For any $t
\in \RR$ set $\mathbf t = (t _ 1, \dots, t _ k)$ with $t _ 1 =
\dots = t _ {k '} = k '  t$ and $ t _ {k ' + 1} = \dots = t _ {k}
= (k ' - k) t$ and $a (t) = \alpha ^ {\mathbf t}$. Then since the
last $k - k '$ entries in the vector $gm$ are zero,
$$
a (t) (gm) \to 0 \qquad \text{as $t \to \infty$},
$$
so that by Mahler's criterion $a (t) x = a (t) g \SL (k, \ZZ) \to \infty$.
\Endproof 
We are finally in a position to finish the proof of Theorem~\ref{theorem about lattice space}:

\vskip6pt {\it Proof of Theorem~{\rm \ref{theorem about lattice space}.}}
Let $H = \prod_ {i = 1 } ^ r \SL (k _ i, \RR)$ be the maximal group fixing
$\mu$, generated by unipotent one-parameter subgroups, and normalized by $A$ as
above. By Theorem~\ref{version of Ratner}, applied to $\mu$ with $\tilde H = A
H$, we know that there is some $L > H$ which is normalized by $A$ so that
almost every $H$-ergodic component of $\mu$ is the $L$-invariant measure on a
closed $L$ orbit. In particular $\mu$ is $L$-invariant, which unless $L<AH$
contradicts the definition of $H$ as the maximal group with the above
properties. Let now $x=g\SL(k,\ZZ)$ have a closed $L$-orbit $Lx$ of finite
volume. Then $\Lambda_L=g^{-1}Lg\cap\SL(k,\ZZ)$ is a lattice in $g^{-1}Lg$, and
so the latter is defined over $\QQ$. Therefore, the same is true for the
semi-simple $g^{-1}Hg=[g^{-1}Lg,g^{-1}Lg]$, $\Lambda_H=g^{-1}Hg\cap\SL(k,\ZZ)$
is a lattice in $g^{-1}Hg$, and $Hx$ is closed with finite volume. However,
this implies $H=L$. 

Thus we conclude that almost every $H$-ergodic component of $\mu$ is supported
on a single $H$-orbit; in other words, in the notations of Lemma~\ref{l: full
rank}, the support of $\mu$ is contained in $X _ H$.

By Lemma~\ref{l: full rank}, this implies that the sum $ \sum_ ik _ i =
k$ since otherwise there is a one-parameter subgroup $a (t)$ of
$A$ so that for every $x \in X_H$ its trajectory $a (t) x \to
\infty$ as $t \to \infty$, in contradiction to Poincar\'e
recurrence.

But if $\sum _ ik _ i = k$, the set $S N _ G (H)$ of Theorem~\ref{version of Ratner} satisfies
$$
S N _ G (H) = N_G(H) = AH
$$
and so by this theorem $\mu$ is supported on a single $AH$-orbit. But $\mu$ is also
$A H$-invariant. This show that $\mu$ is algebraic: an $AH$-invariant
probability measure on a single $AH$-orbit. Note that this $AH$-orbit has
finite volume, hence is closed in $X$.
\hfq

\vskip6pt {\it Proof of Corollary~{\rm \ref{corollary about lattice space}}}.
Let $\mu$ be an $ A$-ergodic probability measure on $X$ with
positive entropy. By Theorem~\ref{theorem about lattice space},
$\mu$ is algebraic, i.e.\ there are a subgroup $A < L < G$ and a
point $x = g \SL (k, \ZZ) \in X$ so that $Lx$ is closed and $\mu$
is the $L$-invariant measure on $Lx$.

Since $\mu$ is a probability measure, this implies that $g ^{-1} L
g \cap \SL (k, \ZZ)$ is a lattice in $g ^{-1} L g$, which, in
turn, implies that $g ^{-1} L g$ is defined over $ \QQ $.
Moreover, the fact that $L$ has any lattice implies it is
unimodular, which in view of $A<L<G$ (and since $A$ is the maximal
torus in $G$) implies $L$ is reductive (this can also be seen
directly from the proof of Theorem~\ref{theorem about lattice
space}).

We conclude that $g ^{-1} L g$ is a reductive group  defined over
$\QQ$, and $g ^{-1} A g$ is a maximal torus in this group. By
\cite [Thm. 2.13]{ Prasad-Raghunathan}, there is an $h \in L$ so
that $g ^{-1} h ^{-1} A hg$ is defined over $ \QQ$ and is
$\QQ$-anisotropic. This implies that $A hx$ is closed and of
finite volume (i.e., since $A \cong \RR ^ {k - 1}$, compact), so
that $Lx$ contains a compact $A$ orbit.

By \cite[Thm. 1.3]{ Lindenstrauss-Weiss}, it follows that
(possibly after conjugation by a permutation matrix), $L$ is the
subgroup of $g=(g_{ij}) \in \SL (k, \RR)$ with $g_{ij}=0$ unless
$i$  is congruent to $ j \bmod m$ for some $ 1 \ne m | k$ (by the
Moore ergodicity theorem it is clear that $A$ acts ergodically on
$Lx$, hence the condition in that theorem that $Lx$ contains a
relatively dense $A$ orbit is satisfied), and that $Lx$ is not
compact. Note that if $k$ is prime this implies that $L = \SL (k,
\RR)$.
\hfq
\vglue16pt
\begin{center}
{\bf Part 2.
 Positive entropy and the set of}\\ {\bf 
 exceptions to Littlewood's Conjecture}
\end{center}

\section{Definitions}\label{sec: definitions}

We recall the definition of Hausdorff dimension, box dimension,
topological and metric entropy. In the following let $Y$ be a
metric space with metric $d_Y(\cdot,\cdot)$.

\Subsec{Notions of dimension}\label{ssec: dimension}
For $D \geq 0$ the {\em $D$-dimensional Hausdorff measure} of a set $B \subset Y$
is defined by
$$
\cH ^ D(B)=\lim_{\varepsilon \rightarrow 0} \inf_{\cC_\varepsilon} \sum_{i}(\operatorname{diam}C_i)^ D,
$$
where $\cC_\varepsilon=\{ C_1,C_2,\dots \}$ is
any countable cover of $B$ with sets $C_i$
of diameter $\operatorname{diam}(C_i)$ less than $\varepsilon$.
Clearly, for $D>m$ any set in the Euclidean space $\RR ^ m$
has Hausdorff measure zero.
The {\em Hausdorff dimension} $\dim_H(B)$ is defined by
\begin{equation}\label{e: def-H-dim}
\dim_H(B)=\inf \{ D:\cH ^ D(B)=0 \}=\sup \{ D:\cH ^ D(B)=\infty \}.
\end{equation}

For every $\varepsilon>0$
a set $F \subset B$ is $\varepsilon$-separated if
$d_Y(x,y)\geq \varepsilon$ for every two different $x,y \in F$. Let
$b_\varepsilon(B)$ be the cardinality of the biggest
$\varepsilon$-separated subset of $B$; then the {\em \/{\rm (}\/upper\/{\rm )}\/ box dimension}
({\em upper Minkowski dimension}) is
defined by
\begin{equation}\label{e: def-box-dim}
\dim_{\mathrm{box}}(B)=
\limsup_{\varepsilon \rightarrow 0}\frac{\log
b_\varepsilon(B)}{|\log \varepsilon|}.
\end{equation}
Note that $b_{\varepsilon_2}(B)\geq b_{\varepsilon_1}(B)$ if $\varepsilon_2<\varepsilon_1$.
Therefore, it is sufficient to consider a sequence
$\varepsilon_n$ in \eqref{e: def-box-dim} if
$\log{\varepsilon_{n+1}}/\log \varepsilon_n \rightarrow 1$ for
$n \rightarrow \infty$.

We recall some elementary properties.
First,  Hausdorff dimension and box dimension do not
change when we use instead of the metric $d_Y(\cdot,\cdot)$
a different but Lipschitz equivalent metric $d'_Y(\cdot,\cdot)$.
The Hausdorff dimension of a countable union is given by
\begin{equation}\label{eq: dimension supremum}
\dim_H\Bigl(\bigcup_{i=1}^ \infty B_i\Bigr)=\sup_i \dim_H(B_i).
\end{equation}
(This follows easily from the fact  that
the measure $\cH ^ D$ is subadditive.)
For any $B$ we have
\begin{equation}\label{eq: Hausdorff smaller box}
\dim_H(B)\leq \dim_{\mathrm{box}}(B).
\end{equation}
If $Y=Y_1 \times Y_2$ is nonempty and $$ d_Y( (x_1,x_2),(y_1,y_2)
)=\max(d_{Y_1}(x_1,x_2),d_{Y_2}(y_1,y_2)),$$ then
$\dim_{\mathrm{box}}Y \leq
\dim_{\mathrm{box}}Y_1+\dim_{\mathrm{box}}Y_2$.

\Subsec{Entropy and the variational principle}\label{ssec: entropy}
Let $T$ be an endomorphism of a {\it compact} metric space $Y$. For
$\varepsilon>0$ and a positive integer $N$ we say that a set $E \subset Y$ is
$(N,\varepsilon)$-separated (with respect to $T$) if for any two different $x,y
\in E$ there exists an integer $0 \leq n<N$ with $d(T ^ nx,T ^ ny)\geq \varepsilon$.
Let $s_{N,\varepsilon}(T)$ be the cardinality of the biggest
$(N,\varepsilon)$-separated set; then the {\em topological entropy} of $T$ is
defined by
\begin{equation}\label{e: def-top-ent}
\h_{\mathrm{top}}(T)=\lim_{\varepsilon \rightarrow 0}\limsup_{N \rightarrow
\infty}\frac{1}{N}\log s_{N,\varepsilon}(T) =\sup_{\varepsilon>0}\limsup_{N
\rightarrow \infty}\frac{1}{N}\log s_{N,\varepsilon}(T).
\end{equation}

Let $\mu$ be a $T$-invariant measure on $Y$, and
let $\cP$ be a finite partition of $Y$ into measurable sets. Then
$$
\Hh_\mu(\cP)=-\sum_{P \in\cP}\mu(P)\log \mu(P)
$$
is the {\em entropy of the finite partition} $\cP$. (Here $0 \log 0=0$.)
For two such partitions $\cP$ and $\cQ$
let $\cP \vee\cQ=\{ P \cap Q:P \in\cP,Q \in\cQ \}$ be the common refinement.
The {\em metric entropy} of $T$ with respect to
$\mu$ and $\cP$ is defined by
$$
\h_\mu(T,\cP)
=\lim_{N \rightarrow \infty}
\frac{1}{N}\Hh_\mu\biggl(\bigvee_{i=0}^{N-1}T ^{-i}\cP\biggr)
$$
and the {\em metric entropy} of $T$ with respect to $\mu$ is
\begin{equation}\label{e: def-metric-ent}
\h_\mu(T) = \sup_\cP \h_\mu(T,\cP),
\end{equation}
where the supremum is taken over all finite partitions $\cP$ of
$Y$ into measurable sets.

Topological and metric entropy are linked:
For a compact metric space $Y$, a continuous map $T:Y \rightarrow
Y$, and a $T$-invariant measure $\mu$ on $Y$ the entropies
satisfy
$$
\h_\mu(T)\leq\htop(T).
$$
Furthermore, the variational principle \cite[Thm.\
8.6]{Walters-82} states that
\begin{equation}\label{e: variational}
\htop(T)=\sup_\mu \h_\mu(T),
\end{equation}
where the supremum is taken over all $T$-invariant measures $\mu$
on $Y$.

\section{Box dimension and topological entropy}\label{sec: partial-hyperbolic}

We return to the study of the left action of the positive diagonal subgroup $A$ on
$X=\SL(n,\RR)/\SL(n,\ZZ)$. We fix an element $a \in A$ and study  
multiplication from the left by $a$ on $X$, in particular we are interested in
the dynamical properties of the restriction $a|_K$ of this map to a compact
subset $K \subset X$. This will lead to a close connection between topological
entropy and box dimension in an {\it unstable manifold}.

The following easy lemma shows that the dimensions can be defined
using the right invariant metric or a norm on $\Mat(k,\RR)$.

\begin{lemma}\label{l: compare-metric}
For every $r>0$ there exists a constant $c_0 \geq 1$ such that
$$
c_0 ^{-1}\| g-h \| \leq d(g,h)\leq c_0 \| g-h \|\mbox{ for all }g,h \in
B_r ^ G,
$$
where $\| A \|=\max_{i,j}|a_{ij}|$ for $A=(a_{ij})\in\Mat(k,\RR)$.
\end{lemma}

$X$ is locally isomorphic to $\SL(k,\RR)$; more specifically, for
every $x \in X$ there exists some $r=r(x)>0$ such that $B_r ^ G$ and
$B_r(x)$ are isomorphic by sending $g$ to $gx$. For small enough
$r$ this is an isometry. For a compact set $K \subset X$ we can
choose $r=r(K)>0$ uniformly with this property for all $x \in K$.

Let $x \in X$, $g \in B_r ^{G}$, $y=gx$, $\bt \in \Sigma$, and
$a=\alpha ^\bt$. Then $ay=(aga ^{-1})ax$. In other words, when we
use the local description of $X$ as above at $x$ and $ax$, left
multiplication by $a$ acts in this local picture like conjugation
by $a$ on $B_r ^ G$. For this reason we define the subgroups
\begin{eqnarray*}
U&=&\{g\in\SL(k,\RR):a^nga^{-n}\rightarrow 0\mbox{ for
}n\rightarrow-\infty\},\\
V&=&\{g\in\SL(k,\RR):a^nga^{-n}\rightarrow 0\mbox{ for
}n\rightarrow\infty\},\mbox{ and}\\
C&=&\{g\in\SL(k,\RR):aga^{-1}=g\},
\end{eqnarray*}
which are the {\em unstable}, {\em stable},  and {\em central} subgroup (for
conjugation with $a$). Let $\delta_{ij}=0$ for $i \neq j$ and $\delta_{ii}=1$,
so that $\Id=(\delta_{ij})$. It is easy to check that $g \in C$ if $g_{ij}=0$
for all $i,j$ with $a_{ii}\neq a_{jj}$, $g \in U$ if $g_{ij}=\delta_{ij}$ for
all $i,j$ with $a_{ii}\geq a_{jj}$, and similarly $g \in V$ if
$g_{ij}=\delta_{ij}$ for all $i,j$ with $a_{ii}\leq a_{jj}$. Furthermore, there
exists a neighborhood $U_0 \subset \SL(k,\RR)$ of the identity so that every $g
\in U_0$ can be written uniquely as $g=g_C g_U g_V$ for some small $g_C \in C$,
$g_{U}\in U$, and $g_V \in V$. If, similarly, $h=h_Ch_Uh_V$, then
\begin{equation}\label{e: metric-bound}
c_1 ^{-1}d(g,h)\leq \max\bigl(d(g_C,h_C),d(g_{U},h_U),d(g_V,h_V)\bigr)\leq
c_1d(g,h)
\end{equation}
for some constant $c_1 \geq 1$.

Since $A$ is commutative, we have $A \subset C$. The map $T(x)=ax$ on $X$ is
{\em partially hyperbolic\/}: $T$ is not hyperbolic (since the identity  is not
an isolated point of $C$), but part of the local description has hyperbolic
structure as follows.

\begin{lemma}\label{l: uniform}
Let $K \subset X$ be compact with $aK \subset K$ and let $r=r(K)$ be as above.
There exists $\lambda>1$ and $c_2 >0$ so that for any small enough
$\varepsilon>0${\rm ,} any $z \in K$ and $f \in B_r ^{U}${\rm ,} and any integer $N \geq 1$
with $d(fz,z)\geq \lambda ^{-N}\varepsilon${\rm ,} there exists a nonnegative integer
$n<N$ with $d(a ^ nfz,a ^ nz)\geq c_2 \varepsilon$.
\end{lemma}

\Proof 
By continuity there exists $\varepsilon \in(0,r)$ such that $d(afa ^{-1},\Id)<r$
whenever $d(f,\Id)<\varepsilon$. This will be the only requirement on $\varepsilon$.
On the other hand, since $U$ is expanded by conjugation with $a$, there exists
some $\lambda>1$ so that $\| afa ^{-1}-\Id \| \geq \lambda \| f-\Id \|$ for all
$f \in B_r ^{U}$. By Lemma~\ref{l: compare-metric} $d(a ^ nfa ^{-n},\Id)\geq
c_0 ^{-2}\lambda ^ n d(f,\Id)$ for all $f \in B_r ^{U}$ and all $n$ for which
$$\max(d(f,\Id ), \dots , d(a ^ nfa ^{-n},\mathrm{Id} ))<r.$$

By assumption $\lambda ^{-N}\varepsilon \leq d(fz,z)=d(f,\Id)<r$.
It follows that there exists $n<N$ with
$c_0 ^{-2}\lambda ^{-1}\varepsilon<d(a ^ nfa ^{-n},\Id)<r$.
Since $a ^ nz \in K$ we get
\vglue8pt
\hfill $
d(a ^ nfz,a ^ nz)=d((a ^ nfa ^{-n})a ^ nz,a ^ nz)=d(a ^ nfa ^{-n},\Id)>c_0 ^{-2}\lambda ^{-1}\varepsilon.
 $
\Endproof\vskip4pt  

We are ready to give a close connection between box dimension and
entropy.

\begin{proposition}\label{prop: Hausdorff step 1}
Let $a \in A$ and $K \subset X$ be compact with $aK \subset K$.
Then one of the following properties holds.
\begin{enumerate}
\item The intersection $Ux \cap K$
of the {\em unstable manifold} $Ux$ with $K$
is a countable union of compact sets of box dimension
zero for every $x \in X$.
\item The restriction $a|_K$ of the multiplication operator
$a$ to $K$ has positive topological entropy.
\end{enumerate}
\end{proposition}

\Proof 
Note, that the first possibility follows
if there exists some $\varepsilon>0$ such
that
\begin{equation}\label{e: claim-C}
P_y=K \cap \bigl(B_{\varepsilon}^{U}y\bigr)
\mbox{ has box dimension zero for every }y \in K.
\end{equation}
To see this, suppose
$K \cap Ux$ for $x \in X$ is nonempty, and cover $K \cap Ux$
by countably many sets $P_y$ as in \eqref{e: claim-C}.
Taking the union for every such $x$ shows the first statement of the proposition.

Now, it suffices to show that if \eqref{e: claim-C} fails for
$\varepsilon$ as in Lemma \ref{l: uniform}, then the topological
entropy $\htop(a|_K)>0$ is positive. Assume $2 \varepsilon \leq r$
and that \eqref{e: claim-C} fails for $y \in K$. We use this to
construct a sequence of $(N,\varepsilon)$-separated sets $F_N \subset
K$. Let $b \in (0,\dim_{\mathrm{box}}(P_y))$. For every $N>0$ let
$F_N \subset P_y$ be a maximal (finite) $\varepsilon
\lambda ^{-N}$-separated set. By choice of $b$ and the definition
of box dimension in \eqref{e: def-box-dim} there are infinitely
many integers $N$ with $|F_N|\geq \lambda ^{bN}\varepsilon ^{-b}$.

We claim that $F_N$ is an $(N,c_2 \varepsilon)$-separated set for $a$
restricted to $K$. Let $gx, hx \in F_N$ be two different points
with $g,h \in B_\varepsilon ^{U}$. By construction $\varepsilon
\lambda ^{-N}\leq d(gx,hx)<2 \varepsilon \leq r$. By Lemma \ref{l:
uniform} applied to $z=hx$ and $f=gh ^{-1}$ there exists a
nonnegative $n< N$ with $d(a ^ ngx,a ^ nhx)\geq c_2 \varepsilon$.
Therefore $F_N$ is $(N,c_2 \varepsilon)$-separated as claimed, and
for infinitely many $N$ we have $s_N(a)\geq|F_N|>\lambda ^{bN}$.
Finally, the definition of topological entropy in \eqref{e:
def-top-ent} implies that $\h_{\mathrm{top}}(\alpha ^\bt)\geq b
\log \lambda>0$.
\Endproof\vskip4pt  

The remainder of this section is only needed for Theorem \ref{theorem about
forms} and Theorem~\ref{thm: transversal-2}. For a compact set which is
invariant in both directions we can also look at the stable and unstable
subgroup simultaneously. Note however, that the set $UV$ is not a subgroup of
$\SL(k,\RR)$.

\begin{lemma}\label{l: uniform 2}
Let $K \subset X$ be compact with $aK=K$. Then
$B ^{U}_{r}B ^{V}_{r}\subset B ^ G_{2r}$ and there exists $\lambda>1$
and $c_3 \geq 1$ so that for any small enough $\varepsilon>0${\rm ,} any $x
\in X$ and $g,h \in B ^{U}_{r}B ^{V}_{r}$ with $hx \in K${\rm ,} and any
integer $N \geq 1$ with $d(gx,hx)\geq \lambda ^{-N}\varepsilon${\rm ,} there
exists an integer $n$ with $d(a ^ ngx,a ^ nhx)\geq c_3 \varepsilon$ and
$|n|<N$.
\end{lemma}

\Proof 
Recall that we use the right invariant metric $d$ to define the balls
$B ^{U}_{r}$, $B ^{V}_{r}$, and $B_{2r}^ G$. Therefore, if
$g_{U}\in B ^{U}_{r}$ and $g_V \in B ^{V}_{r}$, then
$$
d(g_{U}g_V,\Id)\leq
d(g_Ug_V,g_V)+d(g_V,\Id)=d(g_U,\Id)+d(g_V,\Id)<2r
$$
and so $B ^{U}_{r}B ^{V}_{r}\subset B_{2r}^ G$.

If necessary we reduce the size of $r$ such that \eqref{e: metric-bound} holds
for every $g,h \in B_{3r}^ G$. Assume $\varepsilon$ is small enough so that
$aB_{\varepsilon}^ Ga ^{-1}\subset B_{c_1 ^{-1}r}^ G$. Let $\lambda>1$ be such
that $\| afa ^{-1}-\Id \| \geq \lambda \| f-\Id \|$ for $f \in U$ and $\| a
^{-1}fa-\Id \| \geq \lambda \| f-\Id \|$ for $f \in V$.

Let $g,h \in B ^{U}_{r}B ^{V}_{r}$ and $x \in X$ be as in the lemma.
Define $f=gh ^{-1}$, so that $d(f,\Id)\!=\!d(g,h)\geq \lambda ^{-N}\varepsilon$.
Write $f\!=\!f_Cf_Uf_V$ and $w\!=\!\max\bigl(d(f_{U},\Id),d(f_V,\Id)\bigr)$.
By \eqref{e: metric-bound}
\begin{equation}\label{e: estimate-w}
\max\bigl(d(f_C,\Id),w\bigr)\geq c_1 ^{-1}d(f,\Id)
\geq c_1 ^{-1} \lambda ^{-N}\varepsilon.
\end{equation}

We need to rule out the case that $d(f_C,\Id)$ is the only
big term in this maximum. Clearly $f_Ch=f_Ch_Uh_V$ and $g=g_Ug_V$ are the correct
decompositions in the sense of \eqref{e: metric-bound},
and so $d(f_C,\Id)\leq c_1d(f_Ch,g)$. By  right invariance of the metric $d$
we get $d(f_Ch,g)=d(f_C,f)$
and again by \eqref{e: metric-bound} we get that
$d(f_C,f)\leq c_1w$. We conclude
that $d(f_C,\Id)\leq c_1 ^ 2 w$, which allows us to improve \eqref{e:
estimate-w} to $w \geq c_1 ^{-3}\lambda ^{-N}\varepsilon$.

Depending on which term in
$w=\max\bigl(d(f_U,\Id),d(f_V,\Id)\bigr)$ achieves the maximum, we
find either a positive or a negative $n$ with $|n|<N$ so that
$\tilde f=a ^{n}fa ^{-n}=f_C \tilde f_U \tilde f_V$ satisfies
$d(\tilde f, \Id)\in(c_3 \varepsilon ,r)$ for some absolute constant
$c_3$. Since $hx,a ^ nhx \in K$, it follows that
\vglue12pt
\hfill $ d(a ^ ngx,a ^ nhx)=d(\tilde f a ^ nhx,a ^ nhx)=d(\tilde f,\Id)\geq c_3 \varepsilon.$ 
\Endproof\vskip4pt

\begin{lemma}\label{l: Hausdorff_step 2}
Let $a \in A$ and $K \subset X$ be compact with $aK=K$.
Then one of the following properties holds.
\begin{enumerate}
\item
The intersection $B ^{U}_{r}B ^{V}_{r}x \cap K$
has box dimension
zero for every $x \in X$.
\item
The restriction $a|_K$ of the multiplication operator
$a$ to $K$ has positive topological entropy.
\end{enumerate}
\end{lemma}

\Proof 
Suppose that $\dim_{\mathrm{box}}(B_r ^ UB_r ^ Vx \cap K)>b>0$ for
some $x \in X$. By \eqref{e: def-box-dim} there exists a
$\lambda ^{-N}\varepsilon$-separated set $F_N$ for $N \geq 1$, which
satisfies $|F_N|\geq \lambda ^{bN}\varepsilon ^{-b}$ for infinitely
many $N$.

Let  $gx,hx \in F_N$ with $g,h \in H$ and $g \neq h$. By Lemma
\ref{l: uniform 2} there exists an integer $n$ with $|n|<N$ such
that $d(a ^ ngx,a ^ nhx)\geq c_2 \varepsilon$. This shows that
$a ^{-N+1}F_N \subset K$ is $(2N-1,c_2 \varepsilon)$-separated with
respect to $a$. It follows that $s_{2N-1}(a)\geq |F_N|\geq
\lambda ^{bN}\varepsilon ^{-b}$ for infinitely many $N$, and so
$\htop(a|_K)\geq \frac{1}{2}d \log \lambda>0$.
\hfq

\section{Upper semi-continuity of the metric entropy}\label{sec: semi-continuity}
\vglue-6pt

For the construction of an $A$-ergodic measure $\mu$ as in Theorem
\ref{theorem about lattice space} we need one more property of the
metric entropy, namely {\em upper semi-continuity} with respect to
the measure. More specifically we consider the metric entropy
$\h_\mu(a)$ as a function of the $a$-invariant measure $\mu$,
where we use the weak$^*$ topology on the space of probability
measures supported on a fixed compact $a$-invariant set $K$. We
will show that $\limsup_{\ell \rightarrow
\infty}\h_{\mu_\ell}(a)\leq \h_\mu(a)$ whenever $\mu_\ell$ is a
sequence of $a$-invariant measures satisfying $\lim_{\ell
\rightarrow \infty}\mu_\ell=\mu$. This is well known to hold for
expansive maps \cite[Thm.\ 8.2]{Walters-82} and also for
$C ^ \infty$ automorphisms of compact manifolds \cite[Thm.\
4.1]{Newhouse-89}. Strictly speaking, neither of the two results
applies to our case: the left multiplication by $a$ is not
expansive, $X=\SL(k,\RR)/\SL(k,\ZZ)$ is a noncompact manifold, and
there is no reason why the compact subsets $K \subset X$ we study
should be manifolds at all. However, the proof for the expansive
case in \cite[\S8.1]{Walters-82} can be adapted to our
purposes -- which we provide here for the sake of completeness. We
will need a few more facts about entropy and conditional entropy;
see \cite[Ch.~4]{Walters-82} and \cite[Ch.~2 and 4]{Parry-topics}.

Let $\mu$ be a probability measure on a compact metric space $Y$.
Let $\cA \subset\cB_Y$ be a $\sigma$-algebra, which is countably
generated by $A_1,\dots,A_i,\dots\; $. Then the\break atom of $x$ is
defined by \vglue-11pt
$$
[x]_\cA=\bigcap_{i:x \in A_i}A_i \cap \bigcap_{i:x\notin A_i}X \setminus
A_i,
$$
\vglue-2pt
\noindent and the conditional measure $\mu_x ^\cA$ is a probability
measure supported on $[x]_\cA$ a.s. Let $\cP$ be a finite
partition. We will need the notion of conditional entropy \vglue-11pt
$$\Hh_\mu(\cP|\cA)=\int\Hh_{\mu_x ^\cA}(\cP)\operatorname{d}\!\mu$$\vglue-2pt
\noindent
and the following basic properties.

For the trivial $\sigma$-algebra $\mathcal N=\{ \emptyset,Y \}$
the conditional entropy equals the entropy $\Hh_\mu(\cP|\mathcal
N)=\Hh_\mu(\cP)$. For two partitions $\cP$ and $\cQ$ we have the
addition formula
\vglue-14pt
$$
\Hh_\mu(\cP \vee\cQ|\cA)=\Hh_\mu(\cP|\cA)+\Hh_\mu(\cQ|\cP \vee\cA).
$$
\vglue-2pt
\noindent If $\cP$ is finer than $\cQ$ and $\cC\subseteq\cA$ is another
countably generated $\sigma$-algebra, then
$$
\Hh_\mu(\cQ|\cA)\leq\Hh_\mu(\cP|\cC).
$$
Finally, the conditional entropy $\Hh_\mu(\cP|\cA)$ vanishes if
and only if there exists a nullset $N$ such that $[x]_\cA \setminus
N$ is contained in one of the elements of $\cP$ for every $x \in
Y \setminus N$.

Suppose $T:Y \rightarrow Y$ is measure preserving and invertible.
Then the metric entropy \eqref{e: def-metric-ent} of $T$ with
respect to a finite partition $\cQ$ can also be written as
$$
\h_\mu(T,\cQ)=\Hh_\mu\Bigl(\cQ\big|\bigvee_{n=1}^ \infty T ^{-n}\cQ\Bigr);
$$
see \cite[Thm.~4.14]{Walters-82}.\smallbreak

We will also need the dynamical version of relative entropy.
Suppose $\cA$ is a countably generated $\sigma$-algebra that
satisfies $T\cA=\cA$. We define
\begin{equation}\label{e: def-relative-ent}
\h_\mu(T,\cQ|\cA)=\Hh_\mu\Bigl(\cQ\big|\bigvee_{n=1}^ \infty
T ^{-n}\cQ \vee\cA\Bigr);
\end{equation}
then
\begin{equation}\label{e: relative-ent}
\h_\mu(T,\cP \vee\cQ)=\h_\mu(T,\cP)+\h_\mu\Bigl(T,\cQ\big|\bigvee_{i=-\infty}^ \infty T ^{i}\cP\Bigr).
\end{equation}

The entropy with  respect to an invariant measure is defined as a
supremum over all finite partitions; see \eqref{e:
def-metric-ent}. For this reason the following general principle
will be helpful.

\begin{lemma}\label{l: sequence-partition}
Let $a \in A$ and $K \subset X$ be compact with $aK \subset K$. Let
$\mu$ be an $a$-invariant measure supported on $K$.
There exists a sequence of finite partitions $\cQ_m$ of $K$
which satisfies for all $m$
that $\cQ_{m+1}$ is finer than $\cQ_m$. The boundaries
of the elements of $\cQ_m$ are $\mu$-null sets{\rm ,} and the $\sigma$-algebra
$\bigvee_{m=1}^ \infty\cQ_m$ equals the Borel
$\sigma$-algebra $\cB_K$ of $K$. Furthermore{\rm ,}
$\h_\mu(a)=\lim_{m \rightarrow \infty}\h_\mu(a,\cQ_m)$.
\end{lemma}

\Proof 
Let $x \in X$ and define $f(y)=d(x,y)$. Then the measure $f_*\mu$
is a probability measure on $\RR ^+$, there exist arbitrarily small
$\varepsilon>0$ such that $f_*\mu(\{ \varepsilon \})=0$, and so
$\mu(\partial B_\varepsilon(x))=0$.

For $m>0$ we can cover $K$ with finitely many $\varepsilon$-balls
with $\varepsilon<1/m$ whose boundaries are null sets. Let $\cP_m$
be the partition generated by these balls. For $P \in\cP_m$ the
boundary $\partial P$ is contained in the union of the boundaries
of the balls; thus it is a null set. To ensure that the
sequence of partitions is getting finer we define
$\cQ_m=\bigvee_{i=1}^ m\cP_m$. It follows that every $Q \in\cQ_m$ has a
null set as boundary, and that $\cQ_m \nearrow\cB_K$ for
$m \rightarrow \infty$.

The last statement follows from \cite[\S4, Thm.~3]{Parry-topics}.
\hfq

\begin{proposition}\label{p: nice-partition}
Let $a \in A${\rm ,} and $K \subset X$ be compact with $aK \subset K$.
For every finite partition $\cP$ of $K$
into measurable sets with small enough diameters and for
any $a$-invariant measure $\mu$ supported on $K${\rm ,} $\h_\mu(a)=\h_\mu(a,\cP)$.
\end{proposition}

\Proof 
Every $T$-invariant measure $\mu$ which is supported
by $K$ is in fact supported on the set $K'=\bigcap_{n \geq 0}a ^ nK$.
Clearly, $K'$ is compact and satisfies $aK'=K'$.
Since a partition of $K$ into small sets induces
a partition of $K'$ into small sets, we can assume
without loss of generality that $K$ satisfies $aK=K$.

Recall that there exists $r=r(K)>0$ with $d(x,gx)=d(\Id,g)$ whenever
$d(\Id,g)<r$ and $x \in K$.
Let $\varepsilon<r$ be as in Lemma \ref{l: uniform 2}, let $\cP$ be a finite
partition into measurable sets
with diameter less than $\delta$ (to be specified later),
and define the $\sigma$-algebra
$\cA=\bigvee_{i=-\infty}^ \infty a ^{-i}\cP$ generated by the orbit of $\cP$.

Let $x,y \in K$ belong to the same atom of $\cA$; in other words
suppose that for all $i \in \ZZ$ the images $a ^ ix,a ^ iy \in P_i$ belong to the same partition
element of~$\cP$. We claim that (for small enough
$\delta$) this implies $x=f_C y$ for some small $f_C \in C$.

Let $x=fy$ with $f \in B_\delta ^ G$ and suppose $f\notin C$.
Let $f=f_Cf_Uf_V$ with\break $f_C \in C$, $f_U \in U$,
and $f_V \in V$. For small enough $\delta>0$
we have $d(f_C,\Id)<c_3 \varepsilon/2$, $d(f_U,\Id)<r$ and $d(f_V,\Id)<r$.
Let $z=f_Uf_Vy=f_C ^{-1}x$; then $z \in$\break\vskip-12pt\noindent  $B_r ^ UB_r ^ Vy$. Since
$af_C=f_Ca$ we have
$d(a ^ nx,a ^ nz)=d(f_C,\Id)<c_3 \varepsilon/2$ for all $n$.
By Lemma \ref{l: uniform 2} there exists some integer $n$ with
$d(a ^ nz,a ^ ny)\geq c_3 \varepsilon$.  We assume
$\delta<c_3 \varepsilon/2$; then
$$d(a ^ nx,a ^ ny)\geq d(a ^ nz,a ^ ny)-d(a ^ nx,a ^ nz)>c_3 \varepsilon/2 $$ shows
that $a ^ nx$ and $a ^ ny$ cannot belong to the same partition
element of $\cP$. This contradiction shows the claim.

Suppose $\cQ=\{ Q_1,\dots,Q_m \}$ is one of the partitions of Lemma \ref{l:
sequence-partition}. We remove all the boundaries of the elements
of the partition and obtain a partition modulo $\mu$ into open
sets of small diameter.

By \eqref{e: relative-ent} we have
$$
\h_\mu(T,\cQ)\leq \h_\mu(T,\cP \vee\cQ)=\h_\mu(T,\cP)+\h_\mu(T,\cQ|\cA),
$$
where $\h_\mu(T,\cQ|\cA)$ is the relative entropy
as in \eqref{e: def-relative-ent}.
We will show that this last term vanishes, which
together with Lemma \ref{l: sequence-partition}
will conclude the proof of Proposition~9.2.

Let $B \subset X$ be measurable. By Poincar\'e recurrence, there
exists a null set $N$ such that for  every $x\notin N$ and $x \in B$
there exists some $n \geq 1$ with $a ^ n x \in B$. We apply this
simultaneously to the countable family of sets
$$
B_{i,j,g_C,\ell}=\{ x:B_{1/\ell}(x)\subset Q_i \cap g_C ^{-1}
Q_j \}
$$
for $g_C \in C \cap \SL(k,\QQ)$, $Q_i,Q_j \in\cQ$ and $\ell \geq 1$.
To show that the relative entropy
$$
\h_\mu(T,Q|\cA)=\Hh_\mu(\cQ|\tilde\cA)\mbox{ with }
\tilde\cA=\bigvee_{n=1}^ \infty T ^{-n}\cQ \vee\cA (\operatorname{mod} \mu)
$$
vanishes, we have to show that
for $x,y\notin N$ which are in the same atom with respect to $\tilde\cA$
and satisfy $x \in Q_i \in\cQ$ and $y \in Q_j \in\cQ$, in fact, $i=j$
holds. Since $x$ and $y$ belong to the same atom with respect to $\cA$,
we know from the above claim that $y=f_Cx$ for some small $f_C \in C$.
Therefore, $x \in Q_i \cap f_C ^{-1}Q_j$ and there exists some
rational $g_C$ close to $f_C$ with $x \in Q_i \cap g_C ^{-1}Q_j$.
Furthermore, we can ensure that $B_{1/\ell}(x)\subset Q_i \cap
g_C ^{-1}Q_j$, $d(g_C,f_C)<1/\ell$, and $1/\ell<r$.
It follows that $x \in B_{i,j,g_C,\ell}$. By
construction of $N$ there exists $n>0$ with $a ^ nx \in
B_{i,j,g_C,\ell}$. Therefore $a ^ n x \in Q_i$ and
$B_{1/\ell}(a ^ nx)\subset g_C ^{-1}Q_j$. From
$$
d(g_C ^{-1}f_C a ^ nx,a ^ nx)= d(g_C ^{-1}f_C,\Id)<1/\ell<r
$$
we see that $g_C ^{-1}f_C a ^ nx \in g_C ^{-1}Q_j$.
Since $a$ commutes with $f_C$, $f_Ca ^ nx=a ^ ny \in Q_j$.
We have shown that $x \in a ^{-n}Q_i$ and $y \in a ^{-n}Q_j$.
Since $a ^{-n}Q_i,a ^{-n}Q_j$ belong to $\tilde\cA$ and $x,y$ are assumed
to belong to the same atom with respect to $\tilde\cA$, it follows that
$i=j$ as claimed.
\Endproof\vskip4pt  

The above proposition has the following important consequence.

\begin{corollary}\label{c: semi-continuous}
Let $a \in A$ and $K \subset X$ be compact with $aK \subset K$.
Then the metric entropy $\h_\mu(a|_K)$ is
upper semi-continuous with respect to the measure $\mu${\rm {\rm ,}} i.e.{\rm ,} for
every $a$-invariant $\mu$ and every
$\varepsilon>0$ there is a neighborhood $U$ of $\mu$ in the weak$^*$ topology
of probability measures on $K$ such that $\h_\nu(a)\leq \h_\mu(a)+\varepsilon$ for
every $a$-invariant $\nu \in U$.
\end{corollary}

\Proof 
As in the proof of Lemma \ref{l: sequence-partition}
we can find a partition $\cP$ of $K$ whose elements
have small enough diameter to satisfy Proposition \ref{p: nice-partition}
and whose boundaries are null sets with respect to $\mu$.
Therefore $\h_\nu(a)=\h_\nu(a,\cP)$
for every  $a$-invariant measure $\nu$ supported on $K$.
Let $\varepsilon>0$. By the definition
of entropy there exists $N \geq 1$ with
$$
\frac{1}{N}\Hh_\mu\Bigl(\bigvee_{n=0}^{N-1}a ^{-n}\cP\Bigr)<\h_\mu(a,\cP)+\varepsilon/2.
$$
Since the sets in the partition
$\cQ=\bigvee_{n=0}^{N-1}a ^{-n}\cP$ all have boundaries which are
null sets with respect to $\mu$, there exists a weak$^*$
neighborhood $U$ of $\mu$ such that $\nu(Q)$ is very close to
$\mu(Q)$ for every  $Q \in \cQ$. The entropy of the
partition $\cQ$ depends only on the measures of the elements of
$\cQ$; therefore we can make sure that
$$
\frac{1}{N}\bigl|\Hh_\nu(\cQ)-\Hh_\mu(\cQ)\bigr|<\varepsilon/2.
$$
For any $a$-invariant $\nu \in U$,
$$
\h_\nu(a)=\h_\nu(a,\cP)\leq
\frac{1}{N}\Hh_\nu(\cQ)\leq \frac{1}{N}\Hh_\mu(\cQ)+\varepsilon/2<\h_\mu(a,\cP)+\varepsilon,
$$
where we used Proposition \ref{p: nice-partition} for $\nu$ and $\mu$.
Furthermore, $\h_\nu(a,\cP)$ is the infimum over
$\frac{1}{M}\Hh_\nu(\bigvee_{n=0}^{M-1}a ^ n\cP)$ by subadditivity \cite[Thm.~4.10]{Walters-82}.
\hfq

\section{Transversal Hausdorff dimension for the\\ set of points
with bounded orbits}\label{sec: transversal}

In this section we apply Theorem \ref{theorem about lattice space}
to prove two theorems about sets with bounded orbits.

For a unimodular lattice $\Lambda \subset \RR ^ k$ we define
$$\delta_{\RR ^ k}(\Lambda)=\min_{\by \in \Lambda \setminus \{ 0 \}}\| y \|.$$
Clearly, every point $x=m \SL(k,\ZZ)$ with $m \in \SL(k,\RR)$ can be identified
with the unimodular lattice generated by the columns of $m$. By this
identification $\delta_{\RR ^ k}$ becomes a positive continuous function on $X$
with the property that the preimages $K_\rho=\delta_{\RR ^ k}^{-1}[\rho,\infty)$
are compact sets for every $\rho>0$ by Mahler's criterion. In other words $B
\subset X$ is bounded if and only if  $\inf_{x \in B}\delta_{\RR ^ k}(x)>0$.

A nonempty subset $\Sigma'\subset \Sigma$ is a cone if $\Sigma'$ is
convex and satisfies $r\bt \in \Sigma'$ whenever $r>0$ and
$\bt \in \Sigma'$.

\begin{theorem}\label{thm: transversal-1}
Let $X=\SL(k,\RR)/\SL(k,\ZZ)$ with $k \geq 3${\rm ,} and let $\Sigma'$
be an open cone in $\Sigma$.
Define
$$
D=\bigl\{ x \in X:\inf_{\bt \in \Sigma'}\delta_{\RR ^ k}(\alpha ^\bt x)>0\bigr\}
$$
to be the set of points with bounded $\Sigma'$-orbits.
Then for every $\bt \in \Sigma'$ and $x \in X$ the $\alpha ^\bt$-unstable manifold
$Ux$ through $x$ intersects
$D$ in a set $D \cap Ux$ of Hausdorff dimension zero.
In fact{\rm ,} $D \cap Ux$ is a countable union
of sets with upper box dimension zero.
\end{theorem}

\Proof 
For $\rho>0$ we define the compact set
\begin{equation}\label{e: def-d-rho}
D_\rho=\bigl\{ x \in X:\inf_{\bt \in \Sigma'}\delta_{\RR ^ k}(\alpha ^\bt x)\geq \rho\bigr\}.
\end{equation}
Clearly $D=\bigcup_{n=1}^ \infty D_{1/n}$.
Let $\bt \in \Sigma'$, $a=\alpha ^\bt$, and $x \in X$. Then $aD_\rho \subset D_\rho$.
By Proposition \ref{prop: Hausdorff step 1} there are two possibilities;
$D_\rho \cap Ux$ is a countable union of compact sets of box dimension zero, or
$a|_{D_\rho}$ has positive topological entropy.
If the first possibility takes place for all $\rho>0$, the theorem follows from
\eqref{eq: dimension supremum} and \eqref{eq: Hausdorff smaller box}.

We will show that the second possibility cannot happen ever.
Suppose $a|_{D_\rho}$ has positive topological entropy. By the
variational principle (\S \ref{ssec: entropy} and
\cite[Thm.~8.6]{Walters-82}) there exists an $a$-invariant measure
$\nu$ supported on $D_\rho$ with positive metric entropy
$h_\nu(a)>0$. However, we need to find an $A$-ergodic measure with
this property in order to get a contradiction to Theorem
\ref{theorem about lattice space}.

Since $\Sigma'\subseteq \Sigma$ is open we can
find a basis $\bt_1,\dots,\bt_{k-1}\in \Sigma'$ of $\Sigma$.
By construction $K=D_\rho$  is compact and satisfies $\alpha ^\bs K \subset K$ for all
$\bs \in \RR ^+\bt_1+\cdots+\RR ^+\bt_{k-1}$.
For $N>0$ the measure
$$
\nu_N=\frac{1}{N ^{k-1}}\int_0 ^ N \cdots \int_0 ^ N(\alpha ^{s_1\bt_1+\cdots
s_{k-1}\bt_{k-1}})_*\nu\operatorname{d}\!s_1 \cdots\operatorname{d}\!s_{k-1}
$$
is supported on $K$ and $a$-invariant.
Since entropy is affine \cite[Thm.\ 8.1]{Walters-82} and
upper semi-continuous by Corollary \ref{c: semi-continuous} with respect to the
measure,
entropy with respect to a generalized convex combination
of measures is the integral of the
entropies. In particular $\h_{\nu_N}(T)=\h_{\nu}(T)$.

Let $\mu$ be a weak$^*$ limit of a subsequence of $\nu_N$. From the definition
of $\nu_N$ it follows that $\mu$ is $A$-invariant. It is also clear that $\mu$
is supported on $K$. From upper semi-continuity,  the entropy
$\h_\mu(a)\geq \h_\nu(a)>0$ is positive. The ergodic decomposition \eqref{e:
decomposition} of $\mu$ defines $\mu$ as a generalized convex combination of
$A$-ergodic measures $\mu_\tau$, which have, almost surely, support contained in
$K$. Since $\h_\mu(T)>0$, there exists some $A$-ergodic measure $\mu_\tau$ with
$\h_{\mu_\tau}(T)>0$ and support in $K$. This contradicts Theorem~\ref{theorem
about lattice space} and concludes the proof of Theorem~\ref{thm:
transversal-1}.
\Endproof\vskip4pt  

Let $D\subseteq X$ be $A$-invariant. We say $D$ has {\em
transversal box dimension zero} if $\{ g \in B_r ^ G:gx \in D \} $ and
$g_{ii}=1$ for $i=1,\dots ,k$ has box dimension zero for all $x
\in D$. (Note that the particular shape of the set used here does
not matter as long as this set is still transversal to the
subgroup $A$.) It is easy to check, that an $A$-invariant set $D$
with transversal box dimension zero has box dimension $k-1$
(unless $D$ is empty).

\begin{theorem}\label{thm: transversal-2}
Let $X=\SL(k,\RR)/\SL(k,\ZZ)$ with $k \geq 3${\rm ,}
let $A \subset \SL(k,\RR)$ be the subgroup
of positive diagonal matrices.
Define
$$
D=\bigl\{ x \in X:\inf_{a \in A}\delta_{\RR ^ k}(ax)>0\bigr\}
$$
to be the set of points with bounded $A$-orbits.
Then $D$ is a countable union of sets with transversal box dimension zero and has Hausdorff 
dimension $k-1$.
\end{theorem}

Clearly $D$ is $A$-invariant, and nonempty since it contains
every periodic $A$-orbit.

\Proof 
As before we define the $A$-invariant compact sets $D_\rho$ as in
\eqref{e: def-d-rho} with $\Sigma'=\Sigma$. Pick an element
$a=\alpha ^\bt \in A$ with $\bt \in \Sigma$, $t_i \neq t_j$ for $i
\neq j$. Then the corresponding central subgroup equals $C=A$ and
$B_r ^ UB_r ^ V$ is transversal to $A$. Let $x \in X$ and $\rho>0$.

We give some conditions on $r>0$.
Our first restriction is that $B_{3r}(x)$ and $B_{3r}^ G$ should
be isometric. Let $O=B_r ^ A \times B_r ^{U}\times B_r ^{V}$ and use
the metric
$$d_O\bigl((f_C,f_U,f_V),(g_C,g_U,g_V)\bigr)=
\max\bigl(d(f_C,g_C),d(f_U,g_U),d(f_V,g_V)\bigr).$$ Furthermore,
define $\psi:O \rightarrow \SL(k,\RR)$ by $\psi(f_C,f_U,f_V)=
f_Cf_Uf_V$ and assume $\psi$ is invertible and Lipschitz in both
directions (as in \eqref{e: metric-bound}).

Let $P=\{(f_C,f_U,f_V)\in O:\psi(f_C,f_U,f_V)x \in D_\rho \}$.
Since $D_\rho$ is\break $A$-invariant, the set $P'=P \cap(\{ \Id \}
\times B_r ^{U}\times B_r ^{V})$ determines $P=\{(f_C,f_U,f_V)\in
O:(\Id,f_U,f_V)\in P'\}$. Clearly  $\psi(P')x=(B_r ^ UB_r ^ Vx)\cap
D_\rho$. By Lemma \ref{l: Hausdorff_step 2} there are two
possibilities; $P'$ has box dimension zero  or $a$ has positive
topological entropy when restricted to $D_\rho$. As in the proof
of Theorem \ref{thm: transversal-1} the latter contradicts Theorem
\ref{theorem about lattice space}. Therefore $D_\rho$ has
transversal box dimension zero.
\hfq

\section{The set of exceptions to Littlewood's
Conjecture}\label{sec: Littlewood}

For any $u,v \in \RR$, define $\tau _ {u, v}$ to be the point
$$
\tau _ {u, v} =
\begin{pmatrix}
1& 0 & 0\\
u& 1& 0\\
v& 0& 1
\end{pmatrix} \SL (3, \ZZ);
$$
in other words, $\tau _ {u, v}$ is the point in $X$ corresponding
to the lattice in $\RR ^ 3$ generated by $(1, u,v), (0,1,0),$ and
$(0,0,1)$. The following well-known proposition gives the
reduction of Littlewood's conjecture to the dynamical question
which we studied in Section \ref{sec: transversal}; see also
\cite[\S2]{Margulis-97} and \cite[\S30.3]{Starkov-00}.
We include the proof for completeness.

\begin{proposition}\label{p: Littlewood}
The tuple $(u,v)$ satisfies
\begin{equation} \label{Littlewood-2}
\liminf _ {n \to \infty} n \langle {n u} \rangle \langle{n
v}\rangle = 0,
\end{equation}
if and only if the orbit $A ^ + \tau _ {u, v}$ is unbounded where
$A ^ +$ is the semigroup
$$
A ^ + = \left\{ \begin{pmatrix}
e ^{-r-s} &  &  \\
& e ^ r &  \\
&  & e ^ s \\
\end{pmatrix}:\mbox{ for }r,s \in \RR ^+\right\}.
$$
\end{proposition}

\Proof 
By the the properties of $\delta_{\RR ^ k}$
we have to show for $(u,v)\in \RR ^ 2$ that
\eqref{Littlewood-2} holds if and only if
$\inf_{a \in A ^+}\delta_{\RR ^ k}(a \tau_{u,v})=0$.

Suppose
$\varepsilon>0$ and there exists $a \in A ^+$ with
$\delta_{\RR ^ k}(a \tau_{u,v})<\varepsilon$. Then
$$
a \tau_{u,v}=\begin{pmatrix}
e ^{-r-s}& 0 & 0\\
e ^ ru& e ^ r& 0\\
e ^ sv& 0& e ^ s
\end{pmatrix} \SL (3, \ZZ)
$$
and by definition of $\delta_{\RR ^ k}$ there exists nonzero $(n,m_1,m_2)\in \ZZ ^ 3$ with
$$
\Bigl\|\begin{pmatrix}
ne ^{-r-s}\\
ne ^ ru+m_1e ^ r\\
ne ^ sv+m_2e ^ s
\end{pmatrix}\Bigr\|<\varepsilon.
$$
Taking the product of all three entries of this vector we find
that
$$
|ne ^{-r-s}(ne ^ ru+m_1e ^ r)(ne ^ sv+m_2e ^ s)|=|n(nu+m_1)(nv+m_2)|<c \varepsilon ^ 3
$$
is small ($c$ depends only on the norm used in $\RR ^ 3$), and so
\eqref{Littlewood-2} follows. Note that $n \neq 0$ since otherwise the lower two
entries in the vector cannot be small.

Suppose now that  \eqref{Littlewood-2} holds for $(u,v)$.
Let $\varepsilon>0$ and find $n>0$ and $(m_1,m_2)\in \ZZ ^ 2$ with
$|n(nu+m_1)(nv+m_2)|<\varepsilon ^ 5$.
We would like to have additionally that
\begin{equation}\label{e: nice uv}
\max(|nu+m_1|,|nv+m_2|)<\varepsilon.
\end{equation}
Suppose this is not true, and assume without loss of generality that\break
$|nv+m_2|\geq \varepsilon$ and $|n(nu+m_1)|<\varepsilon ^ 4$. Then by Dirichlet's
theorem there exists an integer $q<1/\varepsilon$ so that $\langle qnv
\rangle<\varepsilon$. It follows that $|qn(qnu+qm_1)|<\varepsilon ^ 2$, and
$|qnv+m_2'|<\varepsilon$ for some $m_2'\in \ZZ$. In other words when we replace
$n$ by $nq$ and $m_1,m_2$ by $qm_1$ and $m_2'$ respectively, we see that
\eqref{e: nice uv} and $|n(nu+m_1)(nv+m_2)|<\varepsilon ^ 3$ hold simultaneously.
Therefore we can find $r>0$ and $s>0$ with $e ^ r|nu+m_1|=\varepsilon$ and
$e ^ s|nv+m_2|=\varepsilon$. (If one of the expressions vanishes, we use some large
$r$, resp.\ $s$, instead.) Then $e ^{-r-s}n<\varepsilon$ and $\delta_{\RR ^ k}(a
\tau_{u,v})<c \varepsilon$ follows.
\hfq

\demo{Proof of Theorem {\rm \ref{thm: zero Hausdorff}}}
By Proposition \ref{p: Littlewood} the set $\Bad$ is embedded by the map
$(u,v)\mapsto \tau_{u,v}$ to the set $D$ with $A ^+$-bounded orbits. We apply
Theorem~\ref{thm: transversal-1} with $\Sigma'=\{(-r-s,r,s):r,s>0 \}$.
Therefore $D$ intersects every unstable manifold of $\alpha ^\bt$ in a set of
Hausdorff dimension zero where $\bt=(-2,1,1)$. Note that the unstable manifold
of $\alpha ^\bt$ through $\Id \SL(3,\ZZ)$ is the image of $\tau$. It follows that
$\Bad$ has Hausdorff dimension zero, and similarly that $\Bad$ is a countable
union of sets with box dimension zero.
\hfq

\demo{Proof of Theorem {\rm \ref{theorem about forms}}}
We apply Theorem \ref{thm: transversal-2} and set $\Xi_k=D$. Suppose $m\notin
\Xi_k$; then $\delta_{\RR ^ k}(am)<\varepsilon$ for some $a \in A$. By definition of
$\delta_{\RR ^ k}$ there exists some $\n \in \ZZ ^ k$
such that $\Bigl\|\left(\begin{array}{c}a_{11}m_1(\n)\\ \vdots \\
a_{kk}m_k(\n)\end{array}\right)\Bigr\|<\varepsilon$ and \eqref{e: product of forms}
follows.
\Endproof


\references {910}

\bibitem[1]{Cassels-Swinnerton-Dyer} \name{J.\ W.\ S.\ Cassels} and \name{H.\ P.\ F.\ 
Swinnerton-Dyer},
On the product of three homogeneous linear forms and the
indefinite ternary quadratic forms,
{\it Philos.\ Trans.\ Roy.\ Soc.\ London Ser.\ A\/} {\bf 248} (1955), 73--96.

\bibitem[2]{CFS} \name{I. Cornfeld, S. Fomin}, and \name{Y. Sinai},
{\it Ergodic Theory},
Translated from the Russian by A.\ B.\ Sosinskii, Springer-Verlag Inc., 
New York, 1982.

\bibitem[3]{EinsiedlerKatok02}
\name{M.~Einsiedler} and \name{A.~Katok}, Invariant measures on
{$G/\Gamma$} for split
simple {Lie}-groups {$G$},  {\it Comm.\ Pure Appl.\ Math\/}.\ 
{\bf 56} (2003), 1184--1221.

\bibitem[4]{EinsiedlerKatokNonsplit}
\bibline, Rigidity of measures -- the high entropy
case, and non-commuting foliations, {\it  Israel J.\ Math\/}.\ {\bf
148} (2005), 169--238.

\bibitem[5]{Fed69}
\name{H.~Federer}, {\it Geometric Measure Theory}, Die 
Grundlehren der mathematischen
Wissenschaften, Band 153, Springer-Verlag, New York, 1969.

\bibitem[6]{Furstenberg-1967}
\name{H. Furstenberg},
Disjointness in ergodic theory, minimal sets, and a problem in
{D}iophantine approximation,
{\it Math.\ Systems Theory\/} {\bf 1} (1967), 1--49.

\bibitem[7]{Furstenberg-1961}
\bibline, Strict ergodicity and transformation of the torus,
{\it Amer.\ J.\ Math\/}.\ {\bf 83} (1961), 573--601.


\bibitem[8]{HasselblattKatokSurvey}
\name{B.~Hasselblatt} and \name{A.~Katok}, Principal structures,
in {\it Handbook of Dynamical Systems}, Vol.\ 1A, 1--203,
North-Holland, Amsterdam, 2002.



\bibitem[9]{Host-normal-numbers}
\name{B.\ Host},  Nombres normaux, entropie, translations,
{\it Israel J. Math.\/}
\textbf{91} (1995), 419--428.

\bibitem[10]{Hu94} \name{H.~Hu},
Some ergodic properties of commuting diffeomorphisms,
{\it Ergodic Theory Dynam.\ Systems\/} \textbf{13} (1993),  73--100.

\bibitem[11]{Hurewicz}
\name{W.\ Hurewicz},
Ergodic theorem without invariant measure,
{\it  Ann.\ of Math\/}.\ {\bf 45} (1944), 192--206.

\bibitem[12]{Johnson-invariant-measures}
\name{A.\ S.\ A.\ Johnson}, Measures on the circle invariant under
multiplication by a nonlacunary subsemigroup of the integers, {\it
Israel J.\ Math\/}.\ \textbf{77} (1992),  211--240.

\bibitem[13]{KalininKatok99}
\name{B.~Kalinin} and \name{A.~Katok}, Invariant measures for actions of
higher rank
abelian groups, in {\it Smooth Ergodic Theory and its Applications\/} 
(Seattle, {WA},
1999), 593--637, Amer.\ Math.\ Soc., Providence, RI, 2001.

\bibitem[14]{KalininSpatzier02}
\name{B.~Kalinin} and \name{R.\ J.~Spatzier}, Rigidity of the measurable structure for
algebraic actions of higher-rank Abelian groups, {\it  Ergodic Theory
Dynam.\ Systems\/}   {\bf 25} (2005), 175--200.

\bibitem[15]{KatokSpatzierRigidity}
\name{A.~Katok} and \name{R.\ J.~Spatzier},
Differential rigidity of {A}nosov
actions of higher rank abelian groups and algebraic lattice
action, {\it Tr.\ Mat.\ Inst.\  Steklova\/} \textbf{216} (1997); {\it
Din.\ Sist.\ i Smezhnye Vopr\/}., 
292--319.

\bibitem[16]{KatokSpatzierCocycles}
\bibline,
First cohomology of {A}nosov actions of
higher rank abelian groups and applications to rigidity,
{\it Inst.\ Hautes \'Etudes Sci.\ Publ.\ Math\/}.\ {\bf 79} (1994),  131--156.

\bibitem[17]{KatokSpatzier96}
\bibline, Invariant measures for higher-rank
hyperbolic
abelian actions, {\it Ergodic Theory Dynam.\ Systems\/} \textbf{16} 
(1996), 751--778.

\bibitem[18]{Kleinbock-et-al}
\name{D. Kleinbock, N. Shah}, and \name{A. Starkov},
Dynamics of subgroup actions on homogeneous spaces of {L}ie groups and
applications to number theory, in
{\it Handbook of Dynamical Systems\/}, Vol.\ 1A, 813--930, 
North-Holland, Amsterdam, 2002.

\bibitem[19]{LedrappierYoungII85}
\name{F.~Ledrappier} and \name{L.-S. Young}, The metric entropy of
diffeomorphisms,
{I}{I}.\ {R}elations between entropy, exponents and dimension, {\it
Ann.\ of Math\/}.\ 
\textbf{122} (1985),  540--574.

\bibitem[20]{Lindenstrauss03}
\name{E.~Lindenstrauss}, Invariant measures and arithmetic quantum unique
ergodicity, {\it Ann.\ of Math\/}.\ {\bf  163} (2006), 165--219.

\bibitem[21]{Lindenstrauss-Weiss} \name{E. Lindenstrauss} and \name{B.\ Weiss},
On sets invariant under the action of the diagonal group,
{\it Ergodic Theory Dynam.\ Systems\/} {\bf 21}  (2001),  1481--1500.

\bibitem[22]{RLyons}
\name{R. Lyons}, On measures simultaneously $2$- and $3$-invariant,
{\it Israel J. Math\/}.\ {\bf 61} (1988),  219--224.


\bibitem[23]{Margulis-Oppenheim-proof} \name{G. A. Margulis},
Discrete subgroups and ergodic theory,
in {\it Number Theory{\rm ,} Trace Formulas and Discrete Groups\/} 
(Oslo, 1987), 377--398,
Academic Press, Boston, MA, 1989.

\bibitem[24]{Margulis-97}
\bibline, Oppenheim conjecture,
{\it Fields Medallists\/}' {\it Lectures\/}, {\it World
Sci.\ Ser.\ {\rm 20th} Century Math\/}.\ {\bf 5}, 272--327,
World Sci.\ Publishing, River Edge, NJ,
1997.


\bibitem[25]{Margulis-conjectures}
\bibline, Problems and conjectures in rigidity theory, in
{\it Mathematics\/}: {\it Frontiers and Perspectives\/}, 
161--174, Amer.\ Math.\ Soc.,  Providence, RI, 2000.

\bibitem[26]{Margulis-Tomanov} \name{G. A. Margulis} and \name{G. M. Tomanov},
Invariant measures for actions of unipotent groups over
local fields on homogeneous spaces, {\it Invent.\ Math\/}.\ {\bf 116} (1994), 
347--392.

\bibitem[27]{Margulis-Tomanov-2}
\bibline,
 Measure rigidity for almost linear groups and its applications,
 {\it  J. Anal. Math\/}.\ {\bf 69} (1996),  25--54.

\bibitem[28]{Witte-book}
\name{D. W.\  Morris},  Ratner's theorems on unipotent flows, 
{\it  Chicago Lectures in Mathematics Series\/}, Univ.\ of
Chicago Press, Chicago, IL, 2005.

\bibitem[29]{Witte-rigidity}
\name{D. W.\ Morris},  Rigidity of some translations on homogeneous spaces,
{\it Invent.\ Math\/}.\ \textbf{81} (1985),  1--27.

\bibitem[30]{Mozes-epimorphic}
\name{S.\ Mozes},
 Epimorphic subgroups and invariant measures,
{\it  Ergodic Theory Dynam. Systems} {\bf 15} (1995),
1207--1210.

\bibitem[31]{Newhouse-89}
\name{S.~E. Newhouse}, Continuity properties of entropy, {\it Ann. of
Math\/}.\
\textbf{129} (1989),  215--235.

\bibitem[32]{Hee-unpublished}
\name{H. Oh}, Application of the paper of Einsiedler, Katok and Lindenstrauss on the arithmeticity of some discrete subgroups, preprint, 2 pages.

\bibitem[33]{Hee-TIFR}
\name{H. Oh}, On a problem concerning arithmeticity of discrete groups acting
on $\HH \times \dots \times \HH$, in {\it Algebraic Groups and
Arithmetic\/}, 427--440, Tata Inst.\ Fund.\ Res., Mumbai,
2004.

\bibitem[34]{Parry-topics}
\name{W. Parry}, {\it Topics in Ergodic Theory\/}, {\it Cambridge Tracts in
Mathematics\/} {\bf 75}, Cambridge University Press, Cambridge, 1981.

\bibitem[35]{Pollington-Velani-00}
\name{A.~D. Pollington} and \name{S.~L. Velani}, On a problem in simultaneous
{D}iophantine approximation: {L}ittlewood's conjecture, {\it Acta
Math\/}.\
\textbf{185} (2000),  287--306.

\bibitem[36]{Prasad-Raghunathan}
\name{G. Prasad} and \name{M.\ S. Raghunathan},
Cartan subgroups and lattices in semi-simple groups,
{\it Ann.\ of Math\/}.\ {\bf 96} (1972), 296--317.

\bibitem[37]{Raghunathan-book}  \name{M. S. Raghunathan},
\emph{Discrete Subgroups of Lie Groups},
{\it Ergebnisse der Mathematik und ihrer Grenzgebiete\/},
Band 68, Springer-Verlag, New York, 1972.


\bibitem[38]{Ratner-factors} \name{M. Ratner},
Factors of horocycle flows,
{\it Ergodic Theory Dynam.\ Systems\/} {\bf 2}  (1982),\break 465--489.

\bibitem[39]{Ratner-joinings} 38,
Horocycle flows, joinings and rigidity of products,
{\it Ann. of Math\/}.\ {\bf 118} (1983),\break  277--313.

\bibitem[40]{Ratner-survey} \bibline,
Interactions between ergodic theory, Lie groups, and number theory,
{\it Proc.\ Internat.\ Congress of Math\/}.\  Vol.\ 1, 2
(Z\"urich, 1994), 157--182, BirkhŠuser, Basel, 1995.

\bibitem[41]{Ratner-measure-rigidity}
\bibline, On Raghunathan's measure conjecture,
{\it Ann. of Math\/}.\ {\bf 134}  (1991),  545--607.

\bibitem[42]{Ratner91t}
\bibline, Raghunathan's topological conjecture and
distributions of
unipotent flows, {\it Duke Math. J\/}.\ \textbf{63} (1991),235--280.

\bibitem[43]{Ratner-SL2}  \bibline,
Raghunathan's conjectures for ${\rm SL}(2,\RR)$,
{\it Israel J.\ Math\/}.\ {\bf 80} (1992), 1--31.

\bibitem[44]{Rees}
\name{M.~Rees}, Some {$\mathbb R\sp 2$}-anosov flows, 1982.

\bibitem[45]{Rudolph-2-and-3}
\name{D.\ J.~Rudolph}, {$\times 2$} and {$\times 3$} invariant measures and
entropy, {\it Ergodic Theory Dynam.\ Systems\/} \textbf{10} (1990), 
395--406.

\bibitem[46]{Starkov-00}
\name{A.~N. Starkov}, {\it Dynamical Systems on Homogeneous Spaces},
{\it Translations of Mathematical Monographs\/} {\bf 190}, Amer.\ Math.\ Soc.\, Providence, RI,
2000.

\bibitem[47]{Starkov-strict-ergodicity}
\bibline,
 Minimality and strict ergodicity of homogeneous actions,
 {\it  Mat.\ Zametki} {\bf 66} (1999), 293--301.

\bibitem[48]{Tomanov-maxtori}
\name{G. Tomanov}, Actions of maximal tori on homogeneous spaces, in {\it Rigidity in
Dynamics and Geometry\/} (Cambridge, 2000), 407--424, Springer-Verlag,
New York, 2002.

\bibitem[49]{Walters-82}
\name{P. Walters}, {\it An Introduction to Ergodic Theory\/}, {\it Graduate
Texts in Mathematics\/} {\bf 79}, 
 Springer-Verlag, New York, 1982.

 
\Endrefs
\end{document}